\documentclass[preprint,12pt,3p]{amsart}

\usepackage[top=1in, bottom=1.25in, left=1.25in, right=1.25in]{geometry}
\usepackage{amssymb}
\usepackage{graphics}
\usepackage{epsfig}
\usepackage{epstopdf}

\usepackage{amsmath}
\usepackage{subfig}

\newcommand{\sech}{\mathop{\mathrm{sech}}\nolimits}

\usepackage{dcolumn}% Align table columns on decimal point
\usepackage{bm}% bold math
\usepackage[colorlinks = true,
linkcolor = blue,
urlcolor  = blue,
citecolor = blue,
anchorcolor = blue]{hyperref}% add hypertext capabilities
\usepackage{amsmath}
\newcommand{\ffor}{\text{for }}

 \newtheorem{thm}{Theorem}  
\title[Symmetry breaking]{Symmetry breaking bifurcations in the NLS equation with an asymmetric delta potential}
\author[Rahmi Rusin]{Rahmi Rusin}
\author[Robert Marangell]{Robert Marangell}
\author[Hadi Susanto]{Hadi Susanto}
\begin{document}
\maketitle

\begin{abstract}
We consider the NLS equation with a linear double well potential. Symmetry breaking, i.e., the localisation of an order parameter in one of the potential wells that can occur when the system is symmetric, has been studied extensively. However, when the wells are asymmetric, only a few analytical works have been reported. Using double Dirac delta potentials, we study rigorously the effect of such asymmetry on the bifurcation type. We show that the standard pitchfork bifurcation becomes broken and instead a saddle-centre type is obtained. Using a geometrical approach, we also establish the instability of the corresponding solutions along each branch in the bifurcation diagram.
\end{abstract}

%\begin{keyword}
%% keywords here, in the form: keyword \sep keyword
%nonlinear Schr\"odinger equation  \sep symmetry breaking \sep asymmetric delta potential \sep bifurcation
%% MSC codes here, in the form: \MSC code \sep code
%% or \MSC[2008] code \sep code (2000 is the default)
%PACS numbers: 63.20.Pw, 63.20.Ry,
%\end{keyword}
%\pacs{63.20.Pw, %localized modes
%	63.20.Ry, %anharmonic lattice modes
%}

%\end{frontmatter}

%%
%% Start line numbering here if you want
%%
% \linenumbers

%% main text
\section{Introduction}

Symmetry breaking where an order parameter becomes localised in one of symmetric potential wells, appears as a ubiquitous and important phenomenon in a wide range of physical systems, such as in particle physics \cite{part}, Bose-Einstein condensates \cite{albiez2005direct,zibold2010classical}, metamaterials \cite{liu2014spontaneous}, spatiotemporal complexity in lasers \cite{green1990spontaneous}, photorefractive media \cite{kevrekidis2005spontaneous}, biological slime moulds \cite{sawai2000spontaneous}, coupled semiconductor lasers  \cite{heil2001chaos} and in nanolasers \cite{hamel2015spontaneous}. When such a bifurcation occurs, the ground state of the physical system that normally has the same symmetry as the external potential becomes asymmetric with the wave function concentrated in one of the potential wells. 

A commonly studied fundamental model in the class of conservative systems is the nonlinear Schr\"odinger (NLS) equation. It was likely first considered in \cite{davies1979symmetry} as a model for a pair of quantum particles with an isotropic interaction potential, where the ground state was shown to experience a broken rotational symmetry above a certain threshold value of atomic masses. Later works on symmetry breaking in the NLS with double well potentials include among others \cite{mahmud2002bose, marzuola2010long, jackson2004geometric, susanto11,susanto12}. At the bifurcation point, stable asymmetric solutions emerge in a pitchfork type, while the symmetric one that used to be the ground state prior to the bifurcation, becomes unstable. 

While previous works only consider symmetric double well potentials, an interesting result was presented in \cite{theo}, on a systematic methodology, based on a two-mode expansion and numerical simulations, of how an asymmetric double well potential is different from a symmetric one. It was demonstrated that, contrary to the case of symmetric potentials where symmetry breaking follows a pitchfork bifurcation, in asymmetric double wells the bifurcation is of the saddle-centre type. In this paper, we consider the NLS on the real line with an asymmetric double Dirac delta potential and study the effect of the asymmetry in the bifurcation. However, different from \cite{theo}, our present work provides a rigorous analysis on the bifurcation as well as the linear stability of the corresponding solutions using a geometrical approach, following \cite{jackson2004geometric} on the symmetric potential case (see also \cite{jones1988instability,mara1,mara2} for the approach).

Since the system is autonomous except at the defects, we can analyse the existence %and the stability 
of the standing waves using phase plane analysis. We convert the second order differential equation into a pair of first order differential equations with matching conditions at the defects. In the phase plane, the solution which we are looking for will evolve first in the unstable manifold of the origin, and at the first defect it will jump to the transient orbit, and again evolves until the second defect, and then jumps to the stable manifold to flow back to the origin. We also present the analytical solutions that are piecewise continuous functions in terms of hyperbolic secant and Jacobi elliptic function. We analyse their instability using geometric analysis for the solution curve in the phase portrait. %and obtain the condition for the stability of the standing wave. 

The paper is organised as follows. In Section \ref{model}, we present the mathematical model and set up the phase plane framework to search for the standing wave. In Section \ref{sol}, we discuss the geometric analysis for the existence of the nonlinear bound states and show that there is a symmetry breaking of the ground states. Then, the stability of the states obtained are analysed in Section \ref{stab}, where we show the condition for the stability in terms of the threshold value of 'time' for the standing wave evolving between two defects. In Section \ref{num}, we present our numerics to illustrate the results reported previously. Finally we summarize the work in Section \ref{concl}

\section{Mathematical model}
\label{model}

We consider the one dimensional NLS equation
\begin{equation}
i\psi_t(x,t)=-\psi_{xx}(x,t)+\omega \psi(x,t)-|\psi(x,t)|^2\psi(x,t)+V(x)\psi(x,t),
\label{nls}
\end{equation}
where $\psi\in\mathbb{C}$ is a complex-valued function of the real variables $t$ and $x$. The asymmetric double-well potential $V(x)$ is defined as
\begin{equation}
V (x)=-\delta(x+L)-\epsilon \delta(x-L), \quad 0<\epsilon\leq1,
\label{pot}
\end{equation} 
where $L$ is a positive parameter. We consider solutions which decay to 0 as $x \rightarrow \pm \infty$. The system conserves the squared norm $N=\int_{-\infty}^{\infty}|u(x,t)|^2 \,dx$ which is known as the optical power in the nonlinear optics context, or the number of atoms in Bose-Einstein condensates.

Standing waves of \eqref{nls} satisfy
\begin{equation} \label{stateq}
u_{xx}-\omega u +u^3-V(x)u=0,
\end{equation}
The stationary equation \eqref{stateq} is equivalent to system
$u_{xx}=\omega u-u^3 $ for $x\neq \pm L$ with matching conditions:
\begin{align}
u(\pm L^+)&=u(\pm L^-),\quad u_x(\pm L^+)-u_x(\pm L^-)=-\tilde{V}_\pm u(\pm L), 
\label{mm}
\end{align}
with $\tilde{V}_-=1$ and $\tilde{V}_+=\epsilon$.

Our aim is to study the ground states of \eqref{nls}, which are localised solutions of \eqref{stateq} and determine their stability. We will apply a dynamical systems approach by analysing the solutions in the phase plane. However, before proceeding with the nonlinear bound states, we will present the linear states of the system in the following section.

\section{Linear states}
\label{linstat}

In the limit $u \rightarrow 0$, Eq.\ \eqref{stateq} is reduced to the linear system \begin{equation} \label{linsys}
u_{xx}-\omega u-V(x)u=0,
\end{equation} which is equivalent to the system $u_{xx}=\omega u$  for $x\neq \pm L$ with the matching conditions \eqref{mm}.

\begin{figure}[htbp!]
	\centering
	%	\subfloat[]{\includegraphics[scale=0.5]{L_vs_omega_e1}\label{subfig:L_vs_omega_e1}}\subfloat[]
	{\includegraphics[scale=0.7]{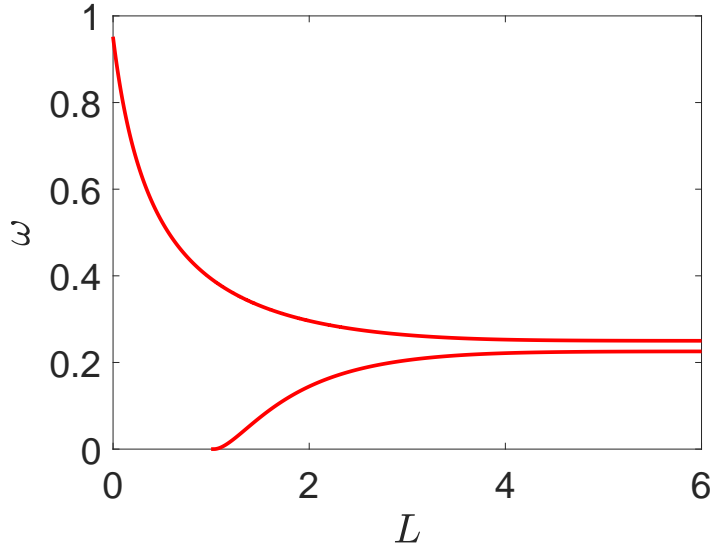}}
	\caption{The eigenvalues $\omega$ as a function of $L$ from \eqref{omega} for $\epsilon=0.95$. The upper curve is $\omega_0$. 	}
	\label{fig:L_vs_omega}
\end{figure}

The general solution of \eqref{linsys} is given by
\begin{equation}
u(x)=
\left\{
\begin{array}{lll}
e^{-\sqrt\omega (x+L)}, &x<-L,\\
Ae^{-\sqrt\omega (x+L)}+Be^{\sqrt\omega (x+L)}, &-L<x<L,\\
Ce^{-\sqrt\omega (x-L)}, &x>L.
\end{array}
\right.
\label{tamb}
\end{equation}
Using the matching conditions, the function \eqref{tamb} will be a solution of the linear system when $A=1-1/2\sqrt{\omega}, B=1/2\sqrt{\omega}, $ and $C=(e^{-2 L \sqrt{\omega }} \left(2 \sqrt{\omega } e^{4 L \sqrt{\omega }}-e^{4 L \sqrt{\omega }}+1\right))/2\sqrt{\omega},$ and  $\omega$ satisfies the transcendental relation
\begin{equation} \label{omega}
L=\frac{1}{4 \sqrt{\omega }}\ln \left(-\frac{\epsilon }{\left(2 \sqrt{\omega }-1\right) \left(\epsilon -2 \sqrt{\omega }\right)}\right).
\end{equation}
This equation determines two bifurcation points of the linear states $\omega_0$ and $\omega_1$. We obtain that the eigenfunction with eigenvalue $\omega_0$ exists for any $L$, while the other one only for $L \geq (1+\epsilon)/2\epsilon$. As $L \rightarrow \infty$, $\omega_0 \rightarrow 1/4$ and $\omega_1 \rightarrow \epsilon^2/4$. We illustrate Eq.\ \eqref{omega} in Fig.\ \ref{fig:L_vs_omega}. Positive solutions that are non-trivial ground state of the system will bifurcate from $\omega_0$, while from $\omega_1$, we should obtain a bifurcation of 'twisted' mode which is not addressed in the present work.

\section{Nonlinear bound states}
\label{sol}

To study nonlinear standing waves (bound states), we convert the second order differential equation \eqref{stateq} into the following first order equations, for $x\neq \pm L,$
\begin{align} \label{ode}
\begin{aligned}
u_x&=y,\\
y_x&=\omega u-u^3, 
\end{aligned}
\end{align}
with the matching conditions
\begin{align}
u(\pm L^+)&=u(\pm L^-),\quad y(\pm L^+)-y(\pm L^-)=-\tilde{V}_\pm u(\pm L) \label{mc1}.
\end{align}

We consider only solutions where $u > 0$. The evolution away from the defects is determined by the autonomous system \eqref{ode} and at each defect there is a jump according to the matching conditions \eqref{mc1}. 

\subsection{Phase plane analysis}
\label{ppa}

System \eqref{ode} has equilibrium solutions $(0,0)$ and $(\sqrt{\omega},0)$ and the trajectories in the phase plane are given by
\begin{equation}
y^2-\omega u^2+\frac{1}{2}u^4=E.
\end{equation}

In the following we will discuss how to obtain bound states of \eqref{nls} which decay at infinity. In the phase plane, a prospective standing wave must begin along the global unstable manifold  $W^u$ of $(0,0)$ because it must decay as $x \rightarrow -\infty$. The unstable manifold $W^u$ is given by \[W^u=\left\{(u,y)|y=\sqrt{\omega u^2-\frac{1}{2}u^4}, 0\leq u \leq \sqrt{2\omega}\right\}.\]

The potential \eqref{pot} will imply two defects in the solutions. After some time evolving in the unstable manifold (in the first quadrant), the solution will jump vertically at the first defect at $x=-L$ according to matching condition \eqref{mc1}. For a particular value of $\omega$, the landing curve for the first jump follows \[J(W^u)=\left\{(u,y)|y=\sqrt{\omega u^2-\frac{1}{2}u^4}-u\right\}.\]

At the first defect, the solution will jump from the homoclinic orbit to an inner orbit as the transient orbit. Let the value of $E$ for the orbit be $\hat{E} \in \left(-\frac{1}{2}\omega^2,0\right)$. If we denote the maximum of $u$ of the inner orbit as $a$, then the value for $a$ in this orbit is \[\hat{a}=\sqrt{\omega+\sqrt{\omega^2+2\hat{E}}}.\]

Denote the value of the solution at the first defect as $u_1$, then it satisfies
\begin{equation} \label{u1val}
u_1^2-2u_1\sqrt{\omega u_1^2-\frac{1}{2}u_1^4}=\hat{E},
\end{equation}
which can be re-written as a cubic polynomial in $u_1^2$,
\[(u_1^2)^3+\left(\frac{1}{2}-2 \omega \right) (u_1^2)^2-\hat{E} u_1^2+\frac{\hat{E}^2}{2}=0.\]
Using Cardan's method \cite{nickalls1993new} to solve the polynomial, we obtain $u_1$ as function of $\hat{E}$, i.e.,   %The product of the three roots of the cubic polynomial is $-\frac{\hat{E}^2}{2}<0$, and since the sum of them, $2\omega-\frac{1}{2}$, then 
for $\omega < 1/4$, there is no real solution, while for $\omega>1/4$, %one of the roots must be negative. Therefore, 
there are 2 real valued $u_1$ given by %satisfying \eqref{u1val}, they are
\begin{align} \label{sol_u1}
\begin{aligned}
u_1^{(1)}&=\left(\frac{1}{3} \left(2 \omega -\frac{1}{2}\right)+\frac{2}{3} \sqrt{3 \hat{E}+\left(\frac{1}{2}-2 \omega \right)^2} \cos \theta\right)^{1/2},\\
u_1^{(2)}&=\left(\frac{1}{3} \left(2 \omega -\frac{1}{2}\right)-\frac{2}{3} \sqrt{3 \hat{E}+\left(\frac{1}{2}-2 \omega \right)^2} \sin \left(\frac{\pi}{6}-\theta\right)\right)^{1/2},
\end{aligned}
\end{align}
where $$\theta=\frac{1}{3} \cos ^{-1}\left(\frac{-54 \hat{E}^2+18 \hat{E} (4 \omega -1)+(4 \omega -1)^3}{\left(12 \hat{E}+(1-4 \omega )^2\right)^{3/2}}\right).$$ For a given $\omega>1/4$, the landing curve of the first jump is tangent to the transient orbit, i.e., $u_1^{(1)}=u_1^{(2)}$, for $\hat{E} = \bar{E}_1$, with
\[\bar{E}_1=\frac{1}{27} \left(36 \omega -\sqrt{(12 \omega +1)^3}-1\right).\] After completing the first jump, the solution will then evolve for `time' $2L$ according to system \eqref{ode}. The `time' $2L$ is the length of the independent variable $x$ that is needed for a solution to flow from the first defect until it reaches the second defect and it will satisfy 
\begin{equation} \label{E1val}
2L=\int_{u_1}^{u_2} \frac{1}{\pm \sqrt{\omega u^2-u^4/2+\hat{E}}} \ du,
\end{equation} 
where $2L=L_1+L_2$, with $L_1$ is the time from $x=-L$ to $x=0$ and $L_2$ is the time from $x=0$ to $x=L$. For $\epsilon=1, L_1=L_2=L$.
The result of the integration of the right hand side of \eqref{E1val} will be in terms of the elliptic integral of the first kind. 

When the solution approaches $x=L$, the solution again jumps vertically in the phase plane according to the matching condition \eqref{mc1}. The set of points that jump to the stable manifold is
\begin{equation}
J^{-1}(W^s)=\left\{(u,y)|y=-\left(\sqrt{\omega u^2-\frac{1}{2}u^4}-\epsilon u \right) \right \}.
\end{equation}

Let $u_2$ be the value of the solution at the second defect. The matching condition \eqref{mc1} gives
\begin{equation} \label{u2val}
\epsilon^2 u_2^2-2\epsilon u_2\sqrt{\omega u_2^2-\frac{1}{2}u_2^4}=\hat{E},
\end{equation}
which can also be re-written as a cubic polynomial in $u_2^2$,
\begin{equation}
u_2^6+\left(\frac{\epsilon^2}{2}-2\omega\right)u_2^4-\hat{E}u_2^2+\frac{\hat{E}^2}{2 \epsilon^2}=0.
\end{equation}
Using a similar argument, the solution exists only for $\omega>\epsilon^2/4$, where in that case the solutions are given by 
\begin{align} \label{sol_u2}
\begin{aligned}
u_2^{(1)}&=\left(\frac{1}{3} \left(2 \omega -\frac{\epsilon^2}{2}\right)+\frac{1}{3} \sqrt{12\hat{E}+\left(\epsilon^2-4 \omega \right)^2} \cos \theta\right)^{1/2},\\
u_2^{(2)}&=\left(\frac{1}{3} \left(2 \omega -\frac{\epsilon^2}{2}\right)-\frac{1}{3} \sqrt{12 \hat{E}+\left(\epsilon^2-4 \omega \right)^2} \sin \left(\frac{\pi}{6}-\theta\right)\right)^{1/2},
\end{aligned}
\end{align}
with $$\theta=\frac{1}{3} \cos ^{-1}\left(-\frac{54 \hat{E}^2+18 \hat{E} \left(\epsilon ^4-4 \omega  \epsilon ^2\right)+\epsilon ^2 \left(\epsilon ^2-4 \omega \right)^3}{\epsilon ^2 \left(12 \hat{E}+\left(\epsilon ^2-4 \omega \right)^2\right)^{3/2}}\right).$$

Similar to the case of the first jump, the landing curve of the second jump is tangent to the transient orbit, i.e., $u_2^{(1)}=u_2^{(2)}$, for $\hat{E} = \bar{E}_2$, with
\[\bar{E}_2=\frac{1}{27} \left(36 \epsilon^2 \omega -\sqrt{\epsilon^2(12 \omega +\epsilon^2)^3}-\epsilon^4\right).\]
For a given value of $\bar{E}_1$ and $\bar{E}_2$, they correspond to certain values of $L$, say $\bar{L}_1$ and $\bar{L}_2$. These values will be used in determining the stability of the solution which will be discussed later in Section \ref{stab}. %For the case $\epsilon=1$, it has been discussed in \cite{jackson2004geometric}. 
For fixed $L$ and $\epsilon$, we can obtain $\hat{E}$ upon substitution of \eqref{sol_u1} and \eqref{sol_u2} to \eqref{E1val} as function of $\omega$, and therefore we can obtain positive-valued bound states for varying $\omega$. 

\begin{figure}[htbp!]
	\centering
	\subfloat[]{\includegraphics[scale=0.5]{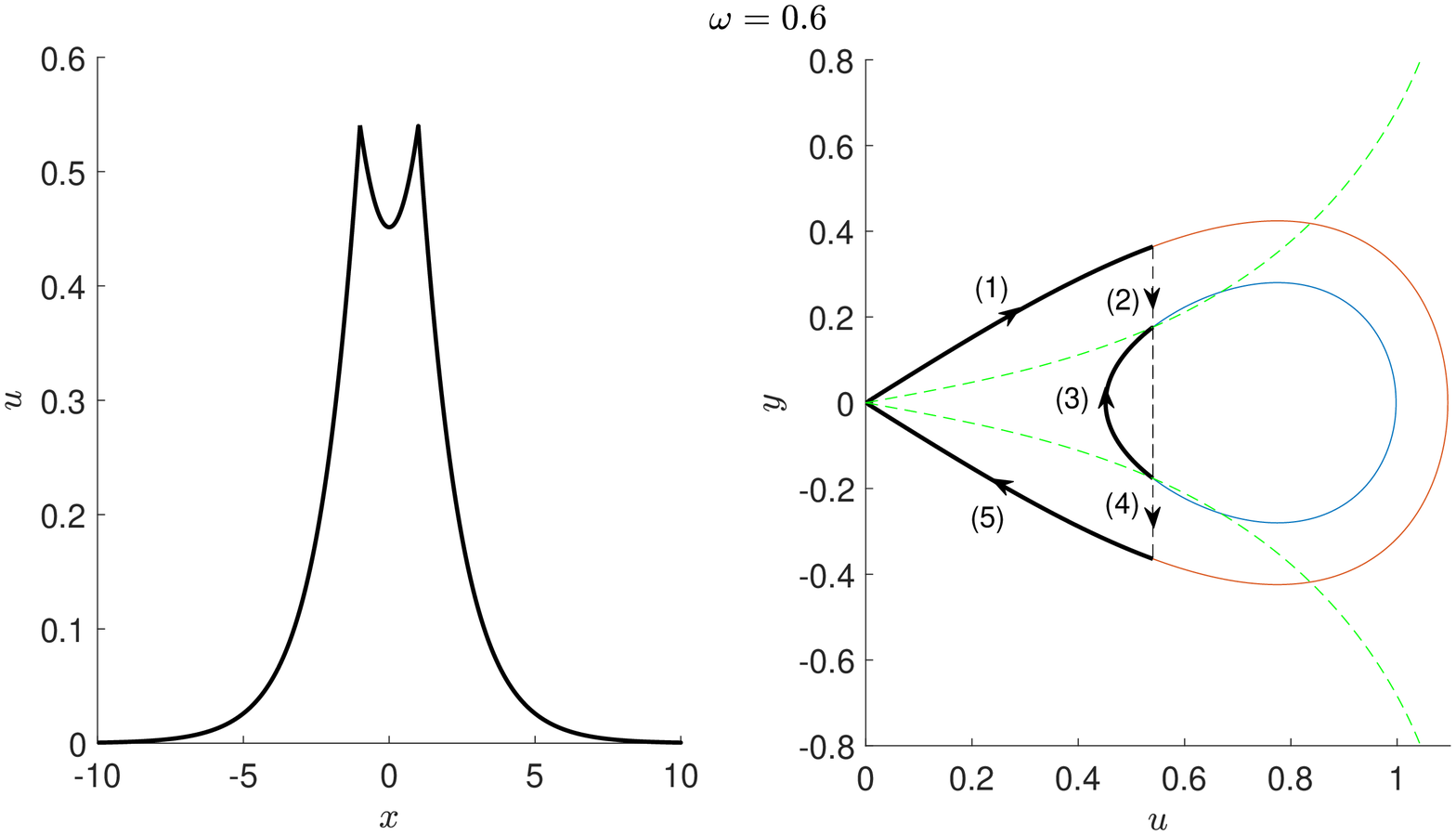}\label{subfig:sym_L_1_e_1_w_0k6}}\\
	\subfloat[]{\includegraphics[scale=0.5]{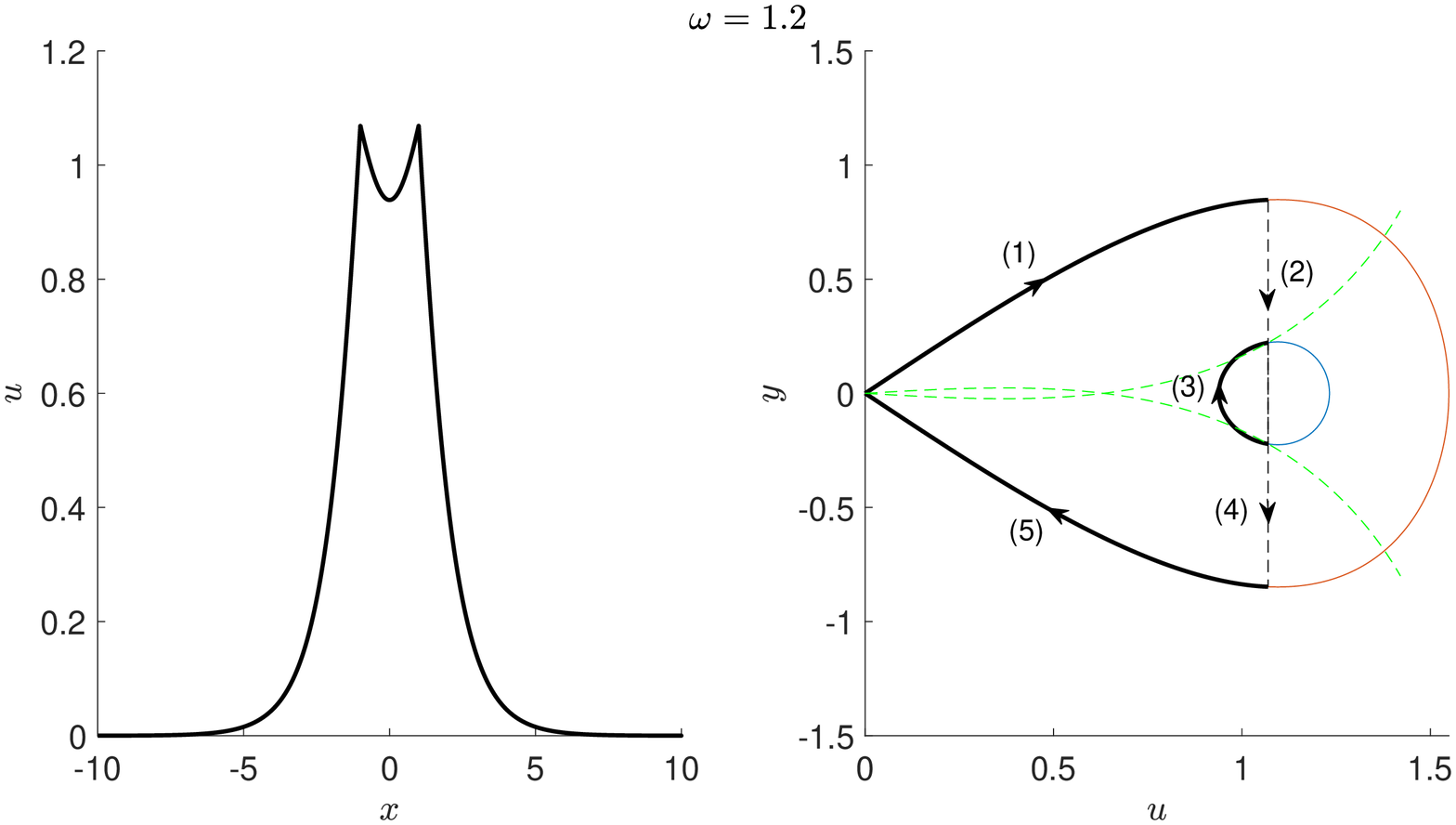}\label{subfig:sym_L_1_e_1_w_1k2}}
	\caption{Localised standing waves of the system for $L=\epsilon=1$ with various values of $\omega$ with $u_1$ and $u_2$ given by (a) $u_1^{(2)}$ and $u_2^{(2)}$,  (b) $u_1^{(1)}$ and $u_2^{(1)}$, (c) $u_1^{(2)}$ and $u_2^{(1)}$, (d) $u_1^{(1)}$ and $u_2^{(2)}$, respectively.
	}
\label{fig:sol}
\end{figure}

\begin{figure}[htbp!]
	\ContinuedFloat
	%\addtocounter{figure}{1}
	\centering
	\subfloat[]{\includegraphics[scale=0.5]{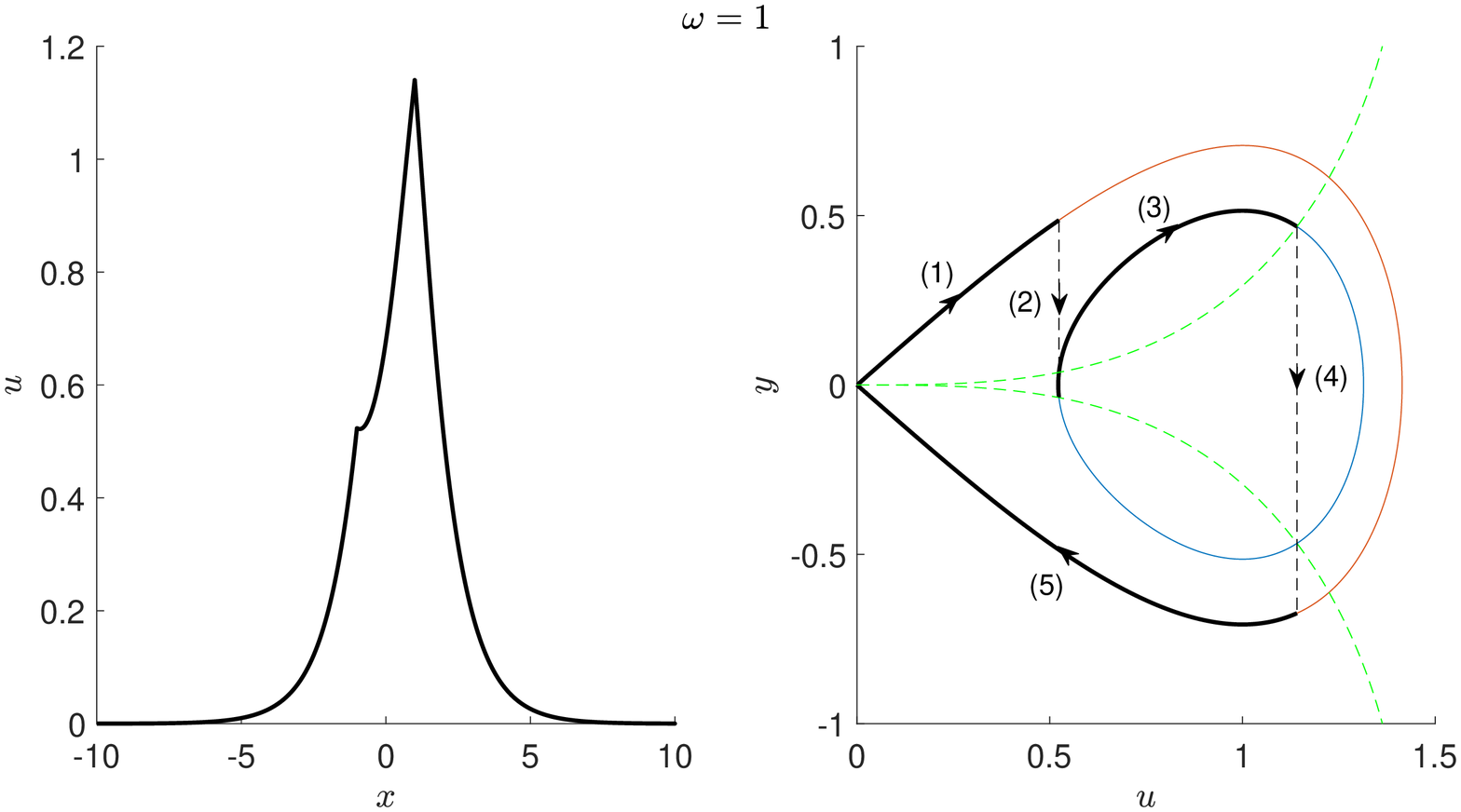}\label{subfig:asym_L_1_e_1_w_1}}\\
	\subfloat[]{\includegraphics[scale=0.5]{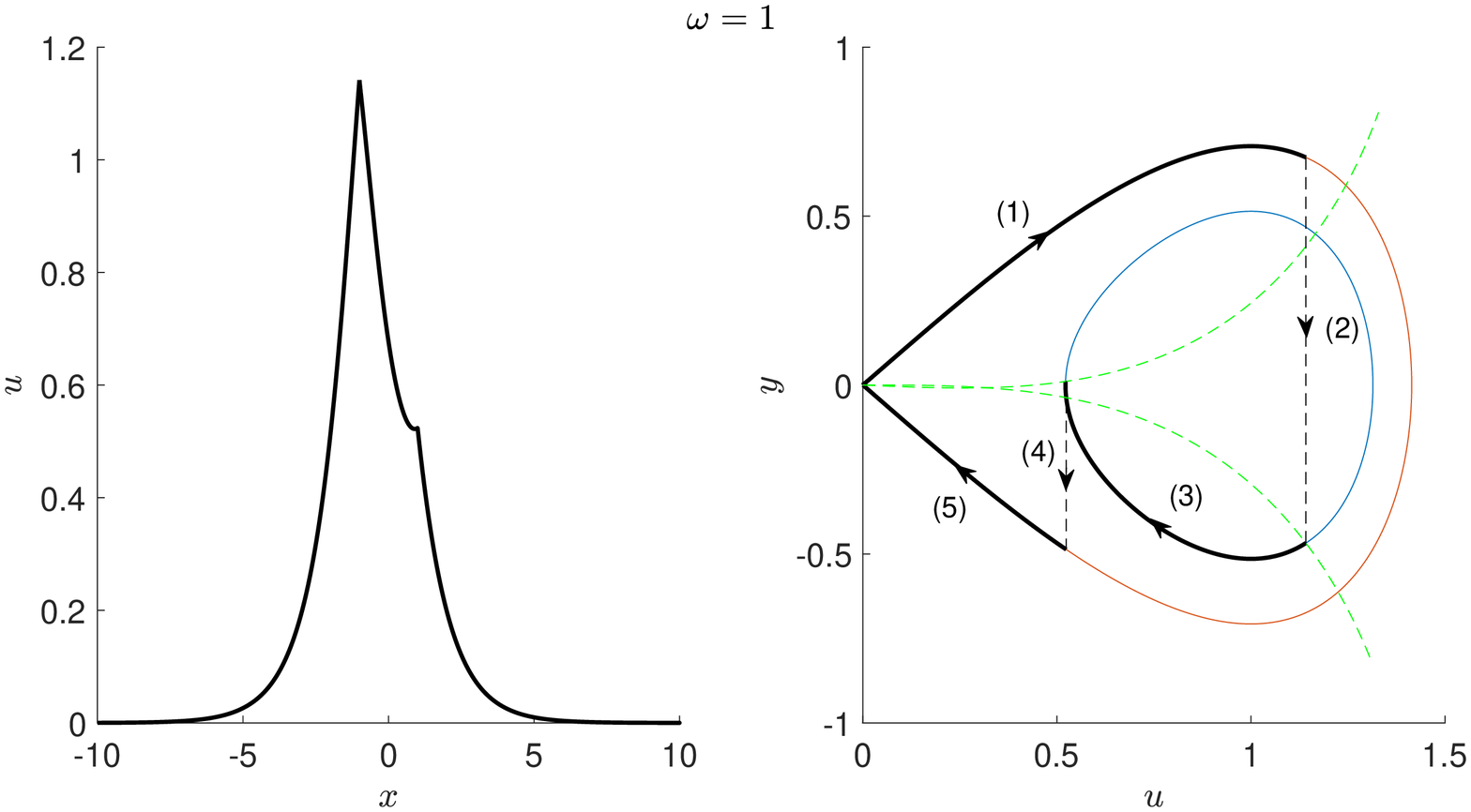}\label{subfig:S2_L_1_e_1_w_1}}
	\caption{(Continued) %Localised standing waves of the system for $L=\epsilon=1$ with various values of $\omega$.
	}
%	\label{fig:solcont}
\end{figure}

\begin{figure}[htbp!]
	\centering
	\subfloat[]{\includegraphics[scale=0.5]{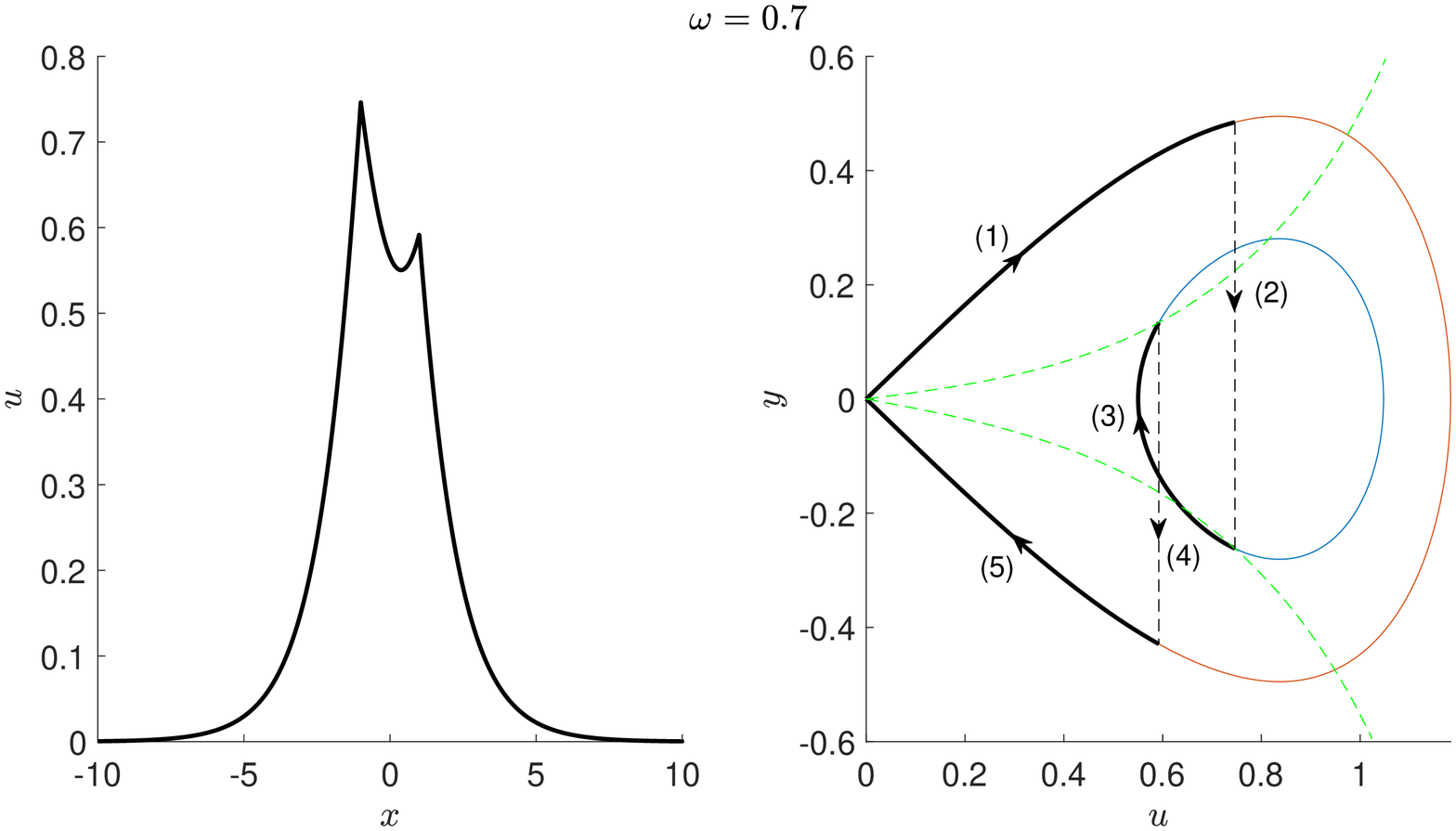}\label{subfig:S2_L_1_e_0k95_w_0k7}}\\
	\subfloat[]{\includegraphics[scale=0.5]{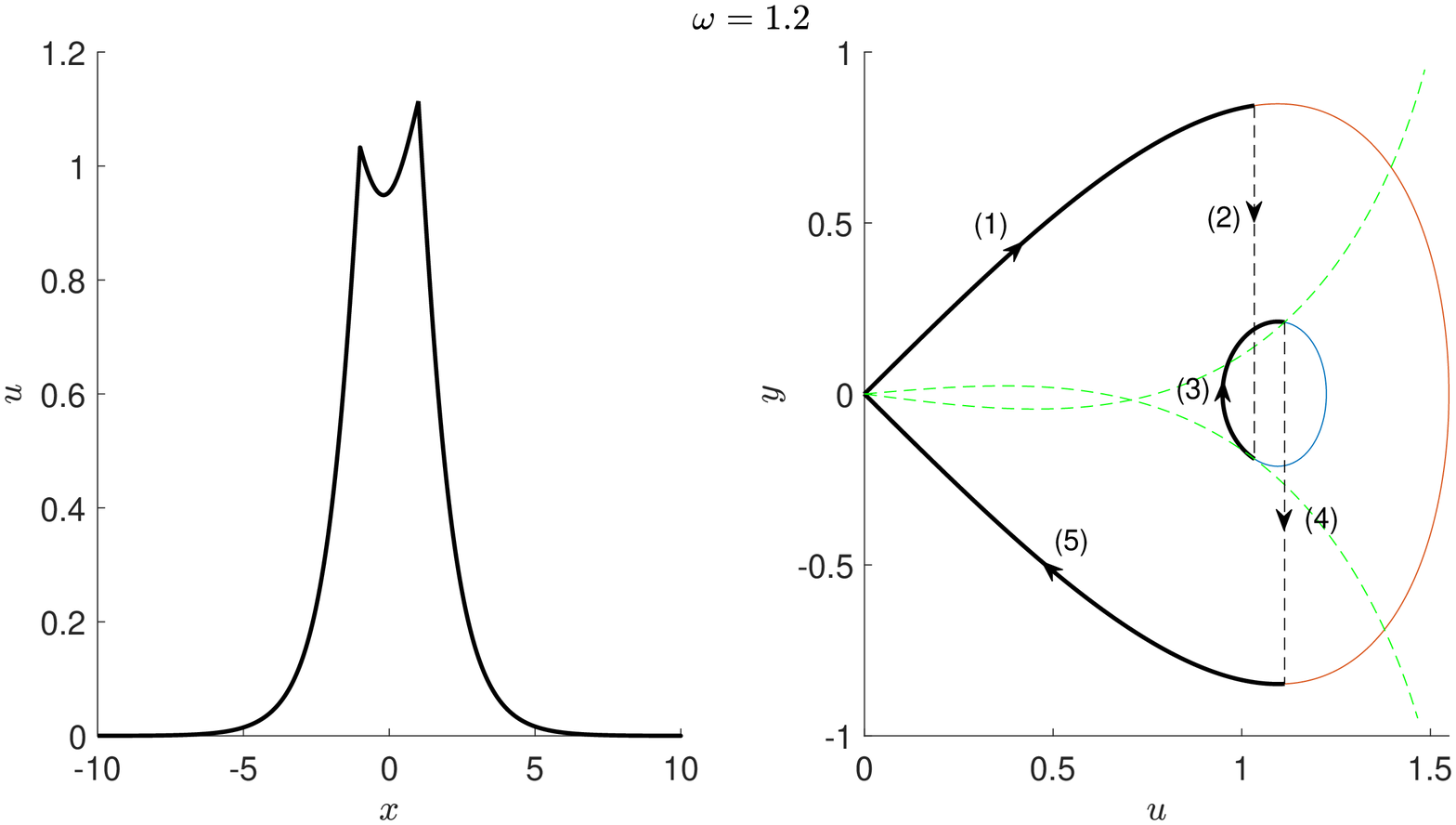}\label{subfig:S1_L_1_e_0k95_w_1k2}}
	\caption{The same as Fig.\ \ref{fig:sol}, but for $\epsilon=0.95$ with $u_1$ and $u_2$ given by (a) $u_1^{(1)}$ and $u_2^{(2)}$,  (b) $u_1^{(1)}$ and $u_2^{(1)}$, (c) $u_1^{(2)}$ and $u_2^{(1)}$, (d) $u_1^{(1)}$ and $u_2^{(2)}$, respectively. 
	}
	\label{fig:sol2}
\end{figure}

\begin{figure}[htbp!]
	\ContinuedFloat
	%\addtocounter{figure}{1}
	\centering
	\subfloat[]{\includegraphics[scale=0.5]{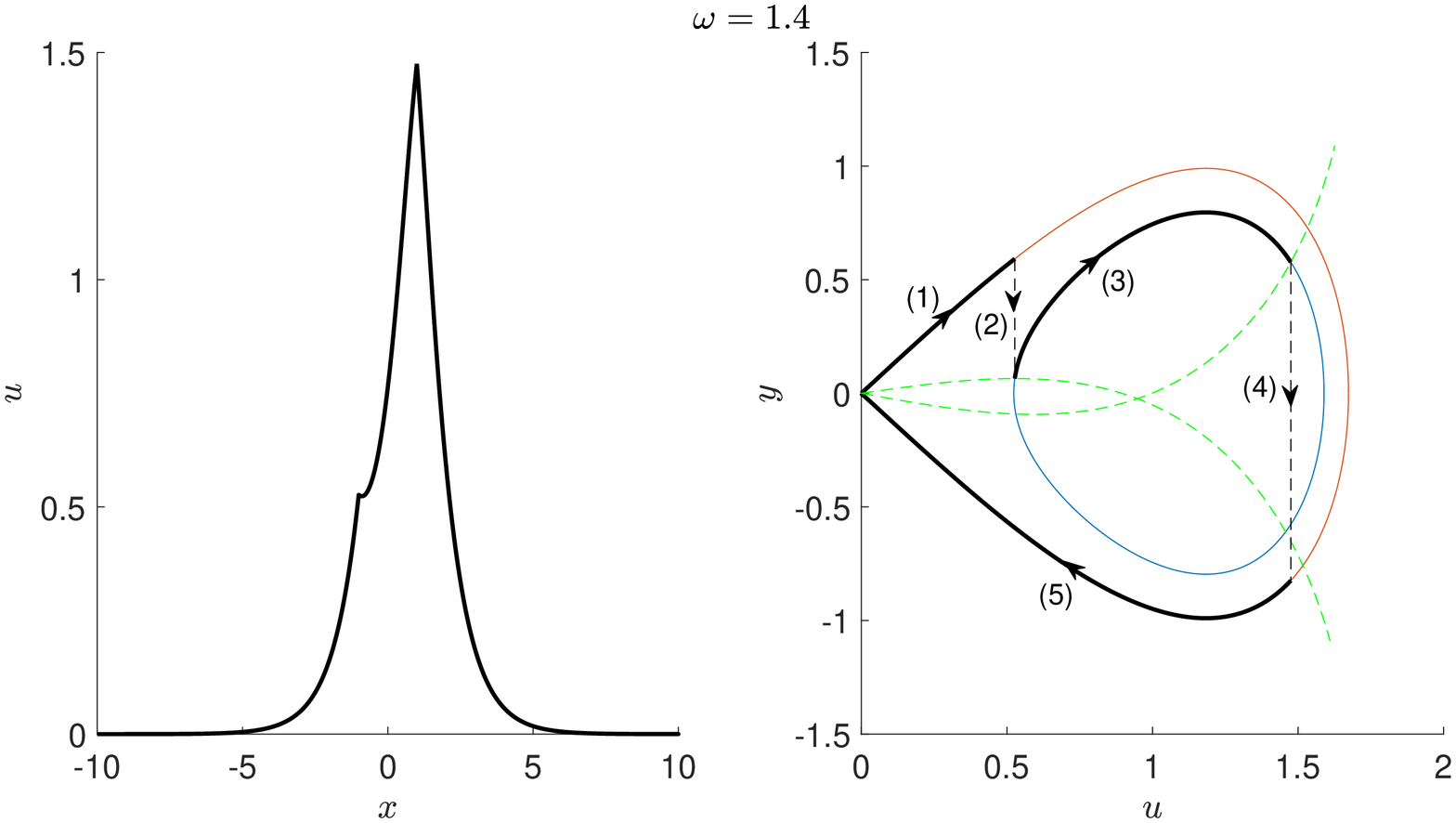}\label{subfig:S3_L_1_e_0k95_w_1k4}}\\
	\subfloat[]{\includegraphics[scale=0.5]{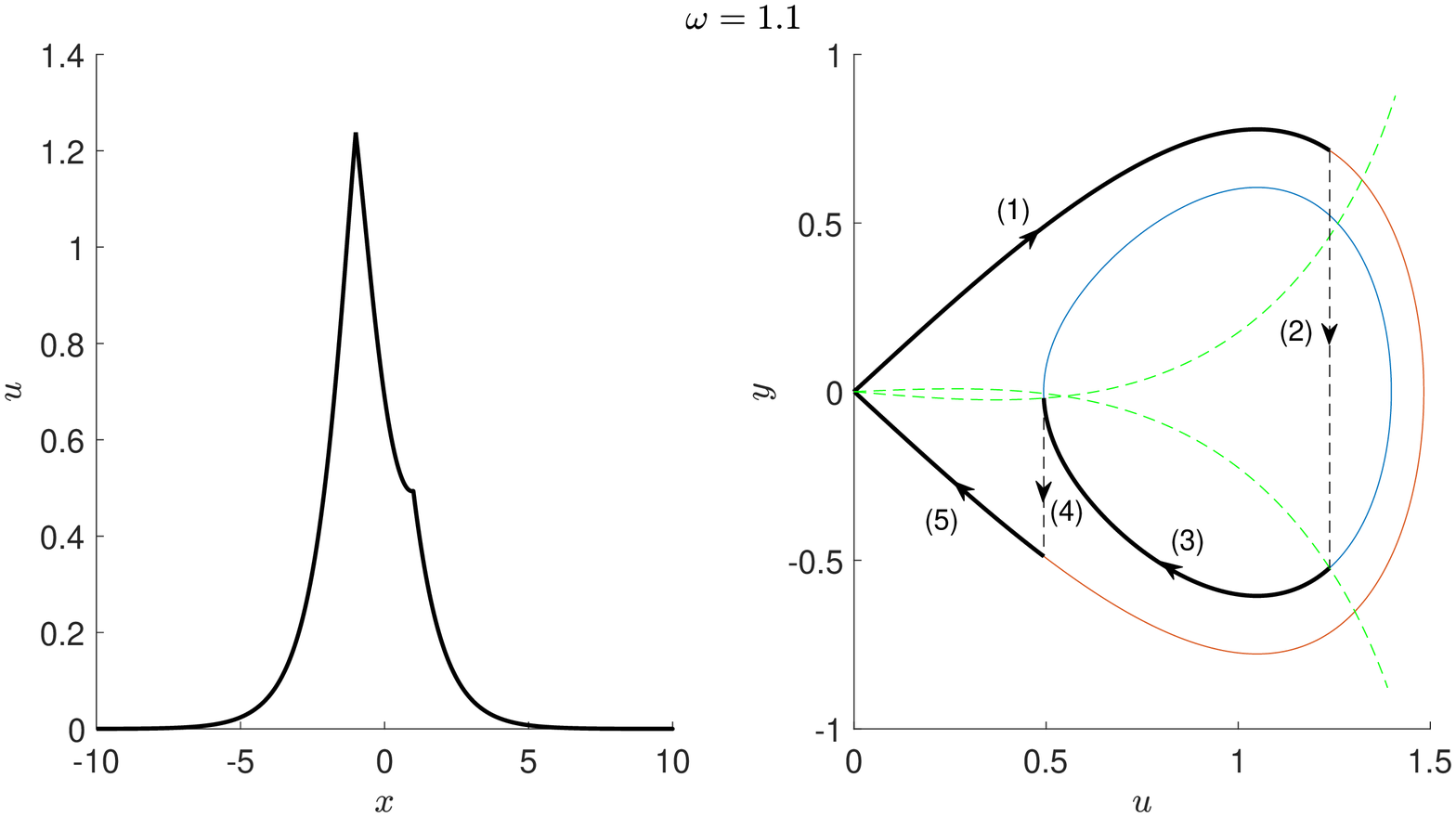}\label{subfig:S2_L_1_e_0k95_w_1k1}}
	\caption{(Continued)
	}
	\label{fig:sol2cont}
\end{figure}

\begin{figure}[tbhp!]
	\centering
	\subfloat[]{\includegraphics[scale=0.5]{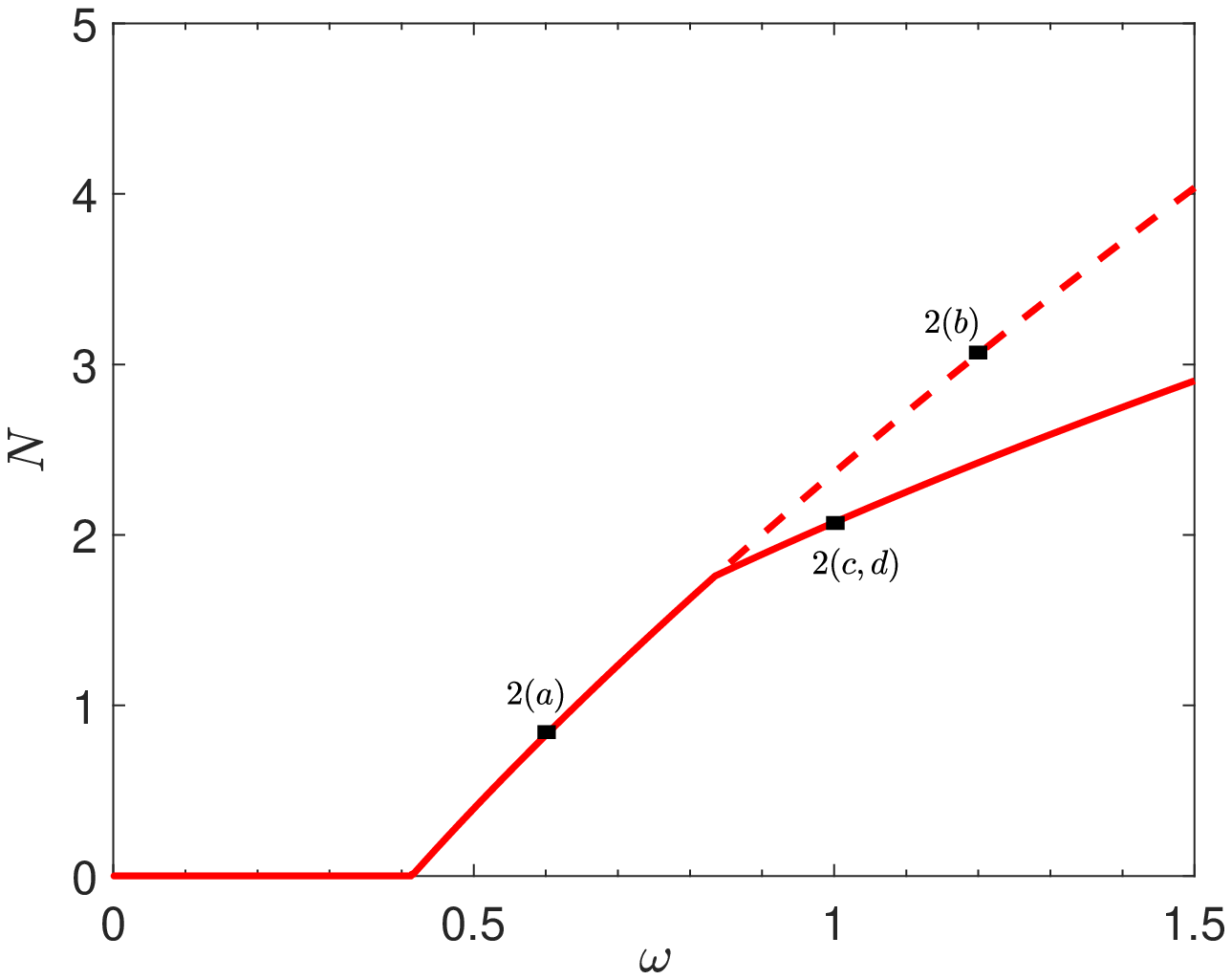}\label{subfig:power_L_1_e_1_num}}
	\subfloat[]{\includegraphics[scale=0.5]{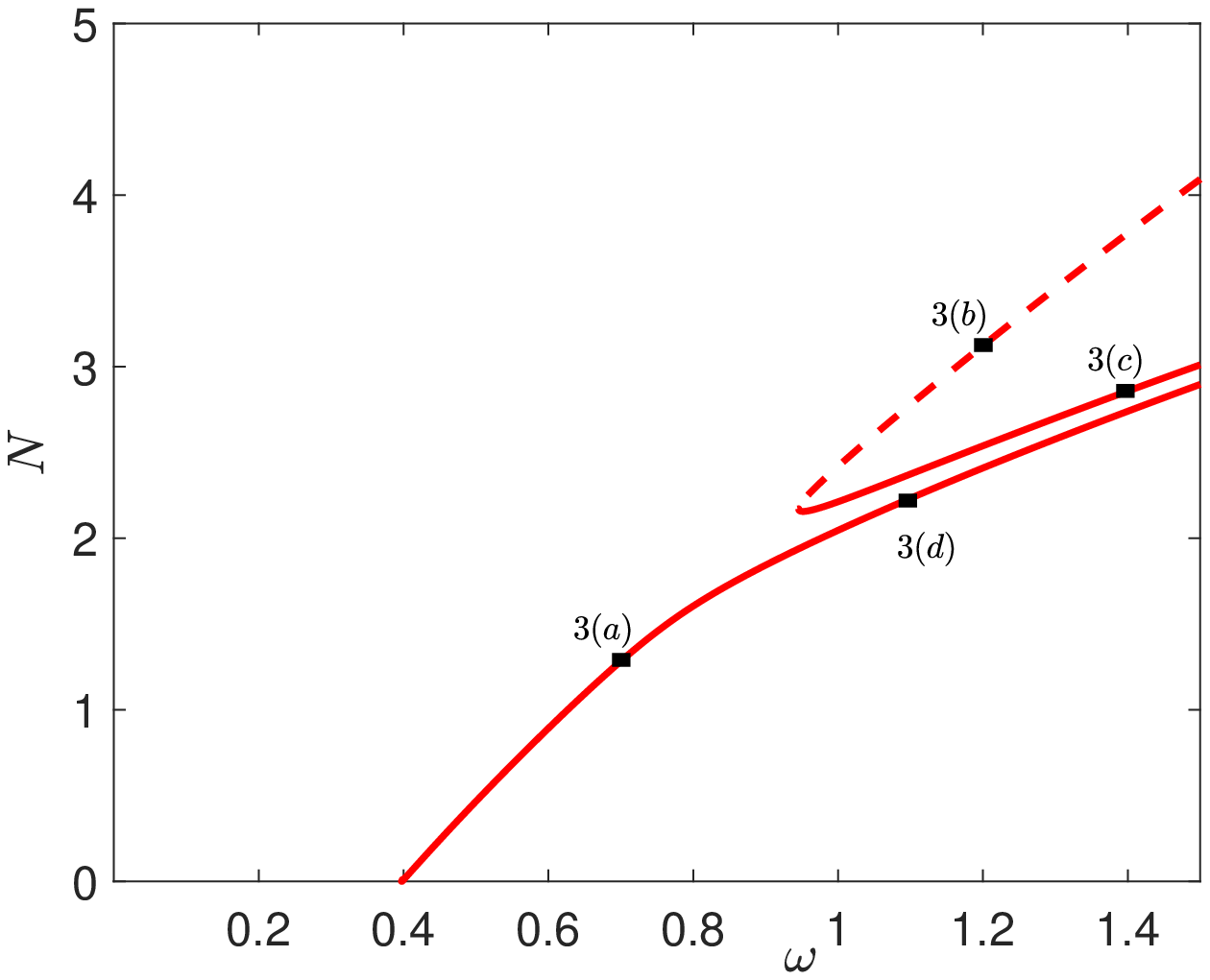}\label{subfig:power_L_1_e_095_num}}
	\caption{Bifurcation diagrams of the standing waves.  Plotted are the squared norms as a function of $\omega$ for $L=1$, and (a) $\epsilon=1$, (b) $\epsilon=0.95$. 
	}
	\label{fig:power}
\end{figure}

We present in Figs.~\ref{fig:sol}--\ref{fig:sol2cont} nonlinear bound states of our system for $L=1$ for $\epsilon=1$ and $\epsilon=0.95$, respectively. We show the solution profiles in the physical space and in the phase plane on the left and right panels, respectively. We also calculate squared norms $N$ of the solutions for varying $\omega$. We plot them in Fig.~\ref{fig:power}. The solid and dashed lines represent the stable and unstable solutions respectively which will be discussed in Section \ref{stab}.

As mentioned at the end of Section \ref{linstat}, indeed standing waves of positive solutions bifurcate from the linear mode $\omega_0$. For $\epsilon=1$, as $\omega$ increases, there is a threshold value of the parameter where a pitchfork bifurcation appears. This is a symmetry breaking bifurcation. Beyond the critical value, we have two types of standing waves, i.e., symmetric and asymmetric states. There are two asymmetric solutions that mirror each other. 

When we consider $\epsilon=0.95$, it is interesting to note that the pitchfork bifurcation becomes broken. The branch of asymmetric solutions splits into two branches and that of symmetric ones breaks into two parts. The upper part of the symmetric branch gets connected to one of the asymmetric branches through a turning point. 

{Using our phase plane analysis, we can determine the critical value of $\omega$ where the bifurcation occurs. The critical value $\omega_c$ as function of $\hat{E}$ can be determined implicitly from the condition when the two roots of $u_1$ \eqref{sol_u1} merge, i.e., 
	\begin{equation}\label{hate}
	\hat{E}=\frac{1}{27} \left(36 \epsilon^2 \omega -\sqrt{\epsilon^2(12 \omega +\epsilon^2)^3}-\epsilon^4\right).
	\end{equation}
Substituting this expression into the integral equation \eqref{E1val}, we can solve it numerically to give us the critical $\omega$ for fixed $L$ and $\epsilon$. %i.e., when $u_1^{(1)}=u_1^{(2)}$ and $u_2^{(1)}=u_2^{(2)}$. 
For $\epsilon=1$, we obtain that $\omega_c \approx 0.8186$ and for $\epsilon=0.95$ we have $\omega_c \approx  0.945$
which agree with the plot in Fig.\ \ref{fig:power}.}

\subsection{Explicit expression of solutions}

The solutions we plotted in Figs.~\ref{fig:sol}  and  \ref{fig:sol2} can also be expressed explicitly as piecewise continuous functions in terms of the Jacobi Elliptic function,  $\text{dn}(rx,k)$. The autonomous system \eqref{ode} has solution \[h(x)=a \, \text{dn}(rx,k)\] with
$r(a,\omega)= \frac{a}{\sqrt{2}} $ and $k(a,\omega)=\frac{2 \left(a^2-\omega\right) }{a^2}$ (for details, see \cite{landa2001regular}). Note that for $a=\sqrt{2\omega}$ we have the homoclinic orbit $h(x)=\sqrt{2\omega}\sech(\sqrt{\omega}x)$. Therefore, an analytical solution of \eqref{stateq} is

\begin{equation}
\label{eq:analyticsol}
u(x)=\begin{cases}
\sqrt{2 \omega}\sech(\sqrt{\omega}(x+\xi_1)), &\ffor \, \, x < -L,\\
a \, \text{dn}(r(x+\xi_2),k), &\ffor \, -L<x<L,\\
\sqrt{2 \omega}\sech(\sqrt{\omega}(x+\xi_3)), &\ffor \, \, x > L,
\end{cases}
\end{equation}
where the constants $\xi_1,\xi_2,$ and $\xi_3$ can be obtained from 
\begin{align}
\begin{aligned}
\xi_1&=\frac{1}{\sqrt{\omega}}\sech^{-1}\left(\frac{u_1}{\sqrt{2\omega}}\right)+L,\\
\xi_2&=\frac{1}{r}\text{dn}^{-1}\left(\frac{u_1}{a},k\right)+L,\\
\xi_3&=\frac{1}{\sqrt{\omega}}\sech^{-1}\left(\frac{u_2}{\sqrt{2\omega}}\right)-L.
\end{aligned}
\end{align}
{The value of $u_1$ and $u_2$ are the same as those discussed in Section \ref{ppa} above. %have to be chosen properly as discussed in  are obtained after $\hat{E}$ from \eqref{E1val} is substituted back to \eqref{sol_u1} and \eqref{sol_u2}. 
	Note that the Jacobi Elliptic function is doubly-periodic. We therefore need to choose the constants $\xi_2$ carefully such that the solution satisfies the boundary conditions \eqref{mm}.}

%\textbf{Tolong dilengkapi biar jelas. Dihubungkan juga dengan solusi yang diplot di Fig 2 dan 3.}

\section{Stability}
\label{stab}

After we obtain standing waves, we will now discuss their stability by solving the corresponding linear eigenvalue problem. We linearise \eqref{nls} about a standing wave solution $\tilde{u}(x)$ that has been obtained previously using the linearisation ansatz $u=\tilde{u}+\delta(p e^{\lambda t}+q^* e^{\lambda^* t})$, with $\delta \ll 1$. %and  the standing wave obtained in Section~\ref{sol}.
Considering terms linear in $\delta$ leads to the eigenvalue problem
\begin{equation}
\lambda\begin{pmatrix}
p \\ q
\end{pmatrix} =\begin{pmatrix}
0 & -\mathcal{L}_- \\
\mathcal{L}_+ & 0
\end{pmatrix} \begin{pmatrix}
p \\ q
\end{pmatrix} =N \begin{pmatrix}
p \\ q
\end{pmatrix}
\label{sys}
\end{equation}
where
\begin{equation}
\begin{aligned}
\mathcal{L}_-&=\frac{d^2}{dx^2}-\omega+\tilde{u}^2-V(x) , \\
\mathcal{L}_+&=\frac{d^2}{dx^2}-\omega +3\tilde{u}^2-V(x).
\end{aligned}
\end{equation}
A solution is unstable when Re$(\lambda)>0$ for some $\lambda$ and is linearly stable otherwise. 

We will use dynamical systems methods and geometric analysis of the phase plane of \eqref{eq:analyticsol} to determine the stability of the standing waves. Let $P$ be the number of positive eigenvalues of $\mathcal{L}_+$ and $Q$ be the number of positive eigenvalues of $\mathcal{L}_-$, then we have the following theorem \cite{jones1988instability}.

\begin{thm}
	If $P-Q \neq 0,1$, there is a real positive eigenvalue of the operator $N$.
	\label{th1}
\end{thm}

The quantities $P$ and $Q$ can be determined by considering solutions of $\mathcal{L}_+p=0$ and $\mathcal{L}_-q=0$, respectively and using Sturm-Liouville theory. The system  $\mathcal{L}_-q=0$ is satisfied by the standing wave $u(x)$, and $Q$ is the number of zeros of standing wave $u(x)$. Since we only consider positive solutions, $Q=0$. By Theorem \ref{th1}, to prove that the standing wave is unstable, we only need to prove that $P\geq 2$. The operator $\mathcal{L}_+$ acts as the variational equation of \eqref{ode}. As such, $P$ is the number of zeros of a solution to the variational equation along $u(x)$ which is `initially' (i.e. at $x=-\infty$) tangent to the orbit of $u(x)$ in the phase plane. Thus $P$ can be interpreted as the number of times initial tangent vector at the origin crosses the vertical as it is evolved under the variational flow. 

Let $\mathbf{p}(u,y)$ be a tangent vector to the outer orbit of the solution at point $(u,y)$ in the phase portrait, and let $\mathbf{q}(u,y)$ be a tangent vector to inner orbit at the point $(u,y)$. That is
\[\mathbf{p}=\begin{pmatrix}
y \\
\omega u-u^3
\end{pmatrix}=\begin{pmatrix}
\pm \sqrt{\omega u^2-\frac{1}{2}u^4} \\
\omega u-u^3
\end{pmatrix},\]
and
\[\mathbf{q}=\begin{pmatrix}
\hat{y} \\
\omega u-u^3
\end{pmatrix}=\begin{pmatrix}
\pm \sqrt{\omega u^2-\frac{1}{2}u^4+\hat{E}} \\
\omega u-u^3
\end{pmatrix}.\]

Let $F$ denote the flow, so $F(\mathbf{p})$ is the image of $\mathbf{p}$ under the flow (together with the matching conditions at the defects). We count the number of times the tangent vector initialised at the origin, say $\mathbf{b}(u,y)$, crosses the vertical as its base point moves along the orbit as $x$ increases. Since the variational flow preserves the orientation of the tangent vector \cite{mara1}, we will use each of the corresponding tangent vectors as the bound of the solution as it evolves after the vector $\mathbf{b}$ is no longer tangent to the orbit due to the defects.  
We will break the orbit into five regions. Let $A_1,A_2,A_3,$ and $A_4$ denote the point $(u_1,y_1)$, $(u_1,y_1-u_1)$, $(u_2, \epsilon u_2+y_2)$, and $(u_2,y_2)$, respectively with $y_1=\sqrt{\omega u_1 -\frac{1}{2}u_1^4}$ and $y_2=-\sqrt{\omega u_2 -\frac{1}{2}u_2^4}$. The first region is  for $x<-L$. On the phase plane it starts from the origin until point $A_1$. The second region is when $x=-L$, i.e., when the solution jumps the first time from $A_1$ to $A_2$. The third region is when $-L<x<L$ where the differential equation \eqref{ode} takes the tangent vector from $A_2$ to point $A_3$. The fourth region is when $x=L$, where the solution jumps for the second time, it jumps from $A_3$ to $A_4$, and last region is for $x>L$ where the vector will be brought back to the origin. 

Let $n_i, i=1,2, \dots, 5,$ denote the number of times $\mathbf{b}$ passes through the verticality in the $i$th region. In the following, we will count $n_i$ in each region. The tangent vector solves the variational flow
\begin{align} \label{varflow}
\begin{aligned}
q_{1,x}&=q_2,\\
q_{2,x}&=q_1-3\tilde{u}^2q_1,
\end{aligned}
\end{align}
where $\tilde{u}$ is the stationary solution.

\subsection{when $x<-L$}
At the first region, we will count $n_1$. It is the region when $\mathbf{b}$ starts from the origin and moves along the homoclinic orbit until it reaches the first defect at $A_1$. At this region, the direction of $\mathbf{b}$ at $(u,y)$ is
\[\tan \theta=\frac{\omega u -u^3}{y}=\frac{\omega-u^2}{\sqrt{1-\frac{1}{2}u^2}}.\]

The sign of $\tan \theta$ depends on the sign of $\omega-u^2$. For $u < \sqrt{\omega}$,  $\tan \theta > 0$, and for $u>\sqrt{\omega}$, $\tan \theta <0$. Since $y>0$, $\mathbf{b}$ points up right in the first quadrant of the plane for the first case, and it points down right in the fourth quadrant for the latter. Therefore, for both cases, the angle must be acute,  $0<|\theta|<\frac{\pi}{2}$. In this part, $n_1=0$. {In what follows we will refer to $\theta$ as {\em the angle of $\mathbf{b}$}}.

\subsection{when $x=-L$}
Next, we will count how many times $\mathbf{b}$ passes through the vertical when it jumps from $A_1$ to $A_2$. The vector $\mathbf{b}$ at $A_2$ is \[\mathbf{b}(A_2)=\begin{pmatrix}
y_1\\
\omega u_1-u_1^3-y_1
\end{pmatrix}
\]  and its direction is $\tan \theta_2=\tan \theta_1 - 1$ with $\theta_i$ the direction of $\mathbf{b}$ in region $i$.  This implies that at the first defect, the vector $\mathbf{b}$ jumps through a smaller angle and larger angle for $L_1<\bar{L}_1$ and $L_1>\bar{L}_1$, respectively. After the jump, $\mathbf{b}$ is tangent to the landing curve $J(W^u)$ but no longer tangent to the orbit of the solution. For $L_1=\bar{L}_1$, after the jump $\mathbf{b}$ will be tangent to the transient orbit but in opposite direction. For all cases, $\mathbf{b}$ does not pass through the vertical. So, up to this stage, $P=n_1+n_2=0$.

\subsection{when $-L<x<L$}
Vector $\mathbf{b}(A_2)$ is now flowed by the equation \eqref {ode} to point $\mathbf{b}(A_3)$.
The variational flow \eqref{varflow} will preserve the orientation of vector $\mathbf{b}$ with respect to the tangent vector of the inner orbit, so $\mathbf{q}$ gives a bound for $\mathbf{b}$ as it evolves. After the first jumping, vector $\mathbf{b}$ points towards into the center (or the concave side) of the inner orbit for $L_1<\bar{L}_1$, and it will point down and out the inner orbit for $L_1>\bar{L}_1$. On the other hand, for $L_1=\bar{L}_1$, the landing curve $J(W^{u})$ is tangent to the inner orbit to which $\mathbf{b}$ jumps, so $\mathbf{b}$ is still tangent to the orbit but now pointing backward. Comparing vector $\mathbf{b}$ with the vector $\mathbf{q}(A_3)$, then up to this point $n_3=0$ or 1.

\subsection{when $x=L$}
In this region, we will count how many times $\mathbf{b}$ cross the vertical when it jumps from $A_3$ to $A_4$, i.e., when $\mathbf{b}(A_3)$ is mapped to $\mathbf{b}(A_4)$. In this region, the tangent vector at $A_3$, $\mathbf{q}(A_3)$ will be mapped to  $F(\mathbf{q}(A_3)$ which has smaller angle, and the jump does not give any additional crossing of the verticality, therefore $n_4=0$

\begin{figure}[htbp!]
	\centering
	\subfloat[]{\includegraphics[scale=0.4]{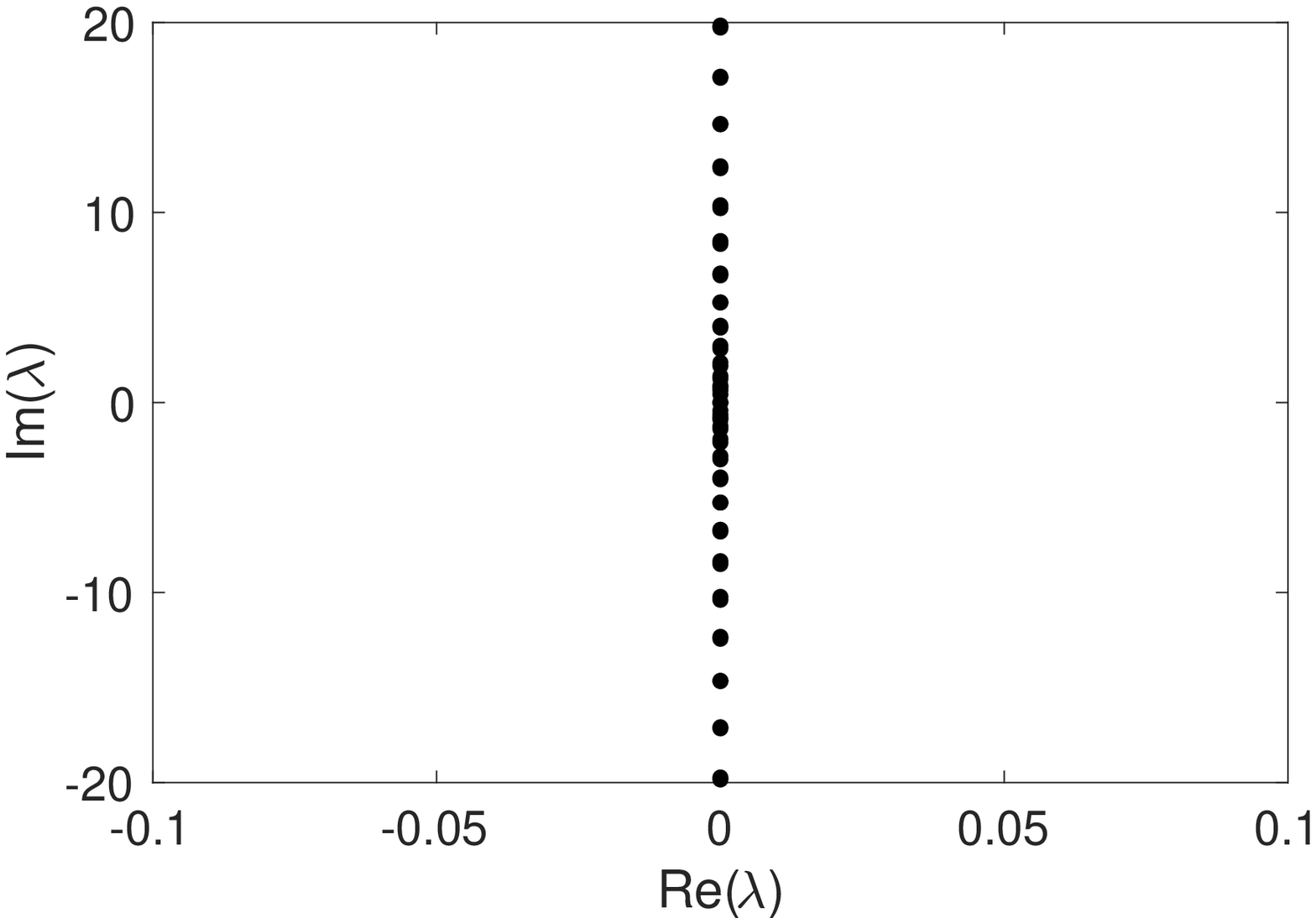}}
	\subfloat[]{\includegraphics[scale=0.4]{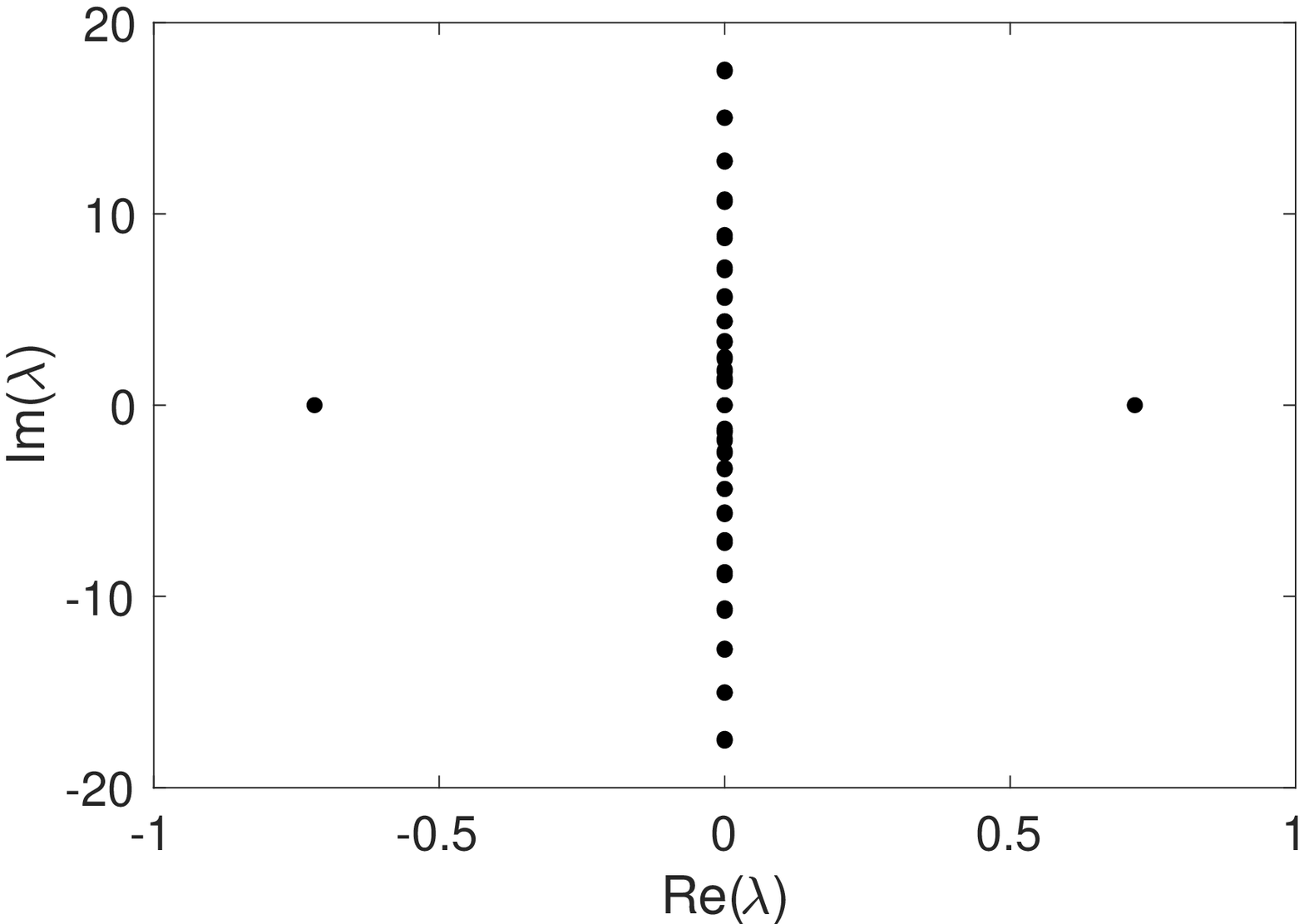}}\\
	\subfloat[]{\includegraphics[scale=0.4]{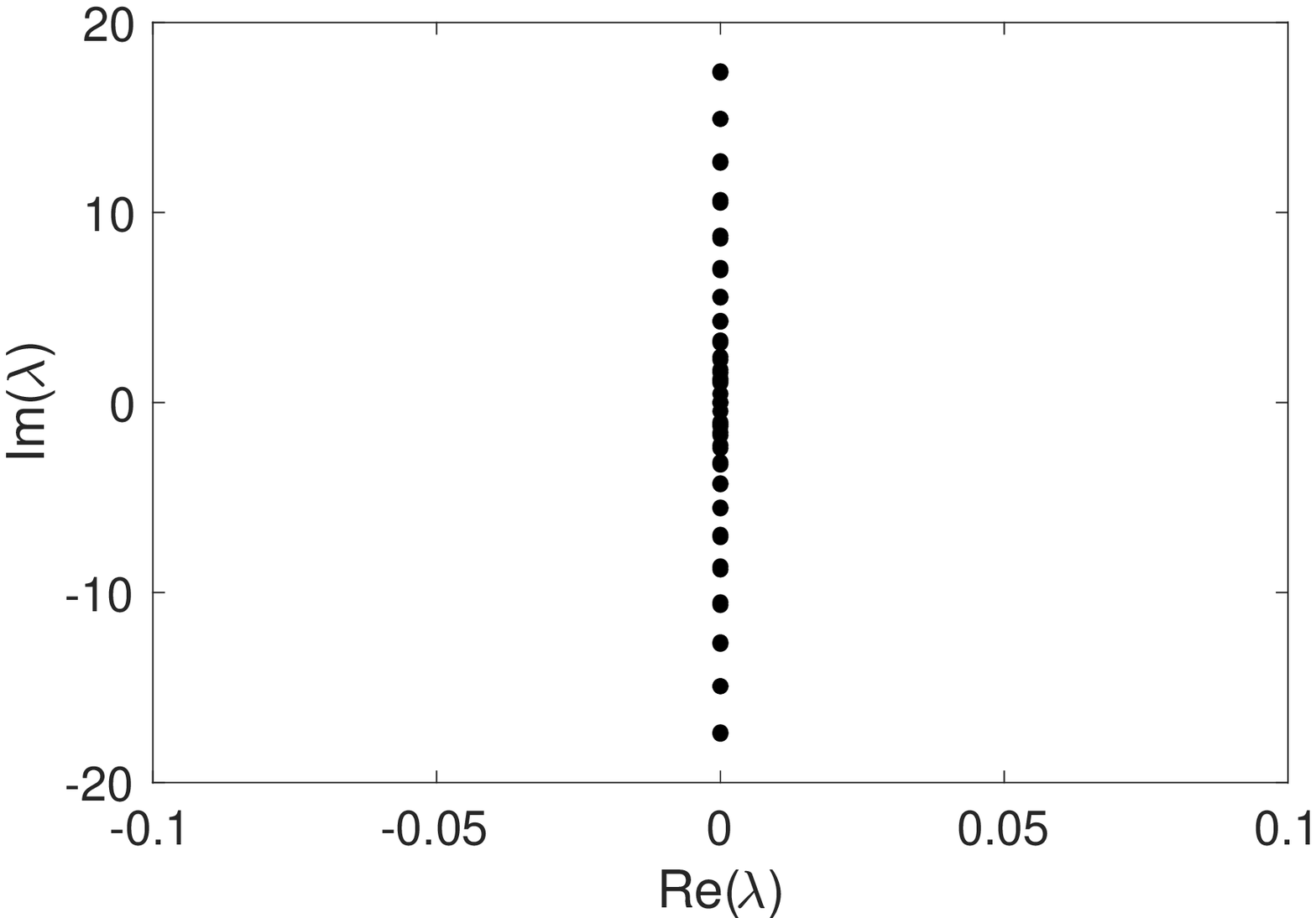}}
	\subfloat[]{\includegraphics[scale=0.4]{as_w1_e1}}
	\caption{Spectrum in the complex plane of the solutions in Fig.\ \ref{fig:sol} in the same order. Panels (c) and (d) %, i.e., the spectrum of Figs.\ \ref{subfig:asym_L_1_e_1_w_1b} and \ref{subfig:asym_L_1_e_1_w_1}, 
		are identical because the solutions are mirror symmetries of each other. 
	}
	\label{fig:sp_sol2}
\end{figure}

\begin{figure}[tbhp!]
	\centering
	\subfloat[]{\includegraphics[scale=0.4]{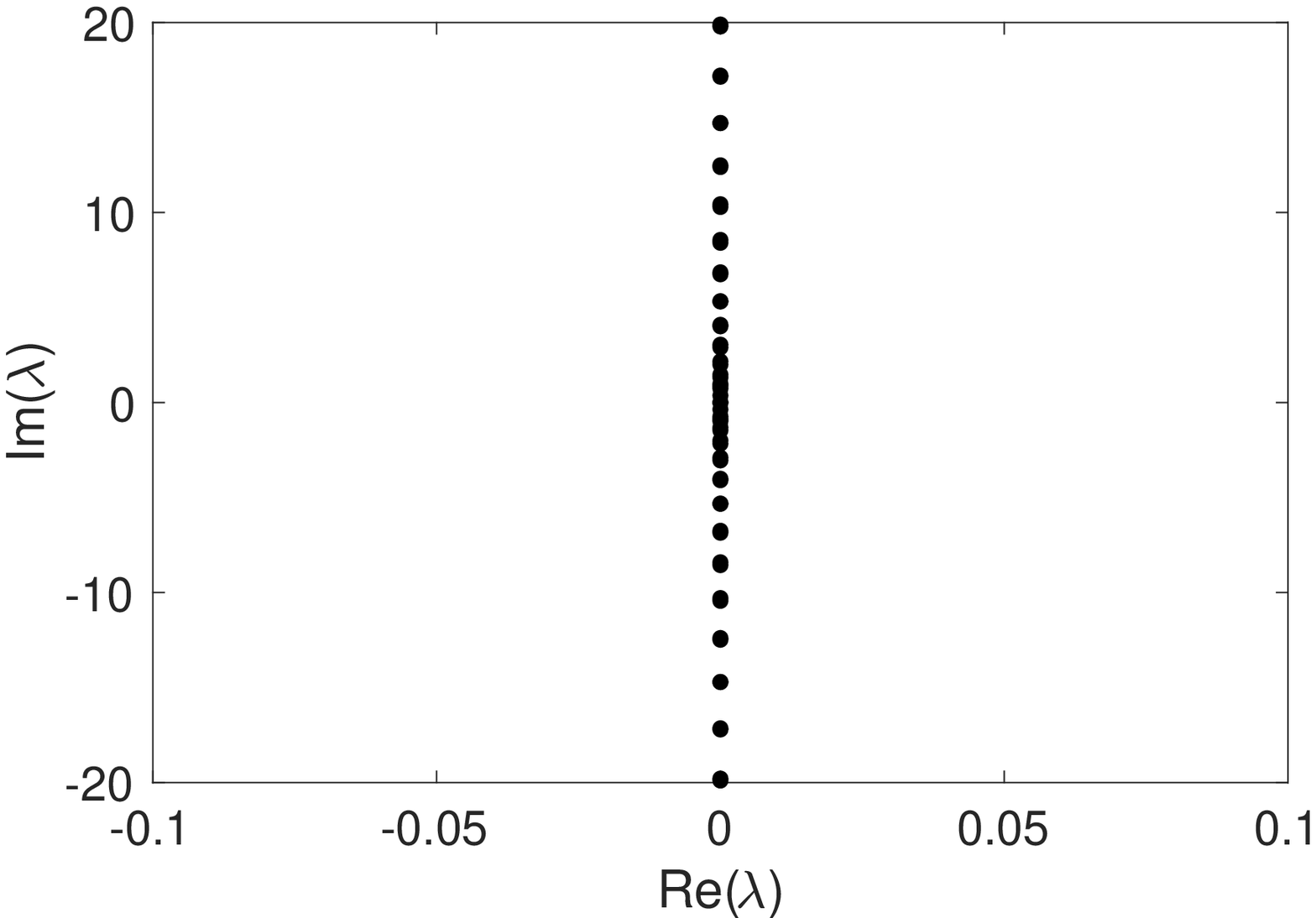}}
	\subfloat[]{\includegraphics[scale=0.4]{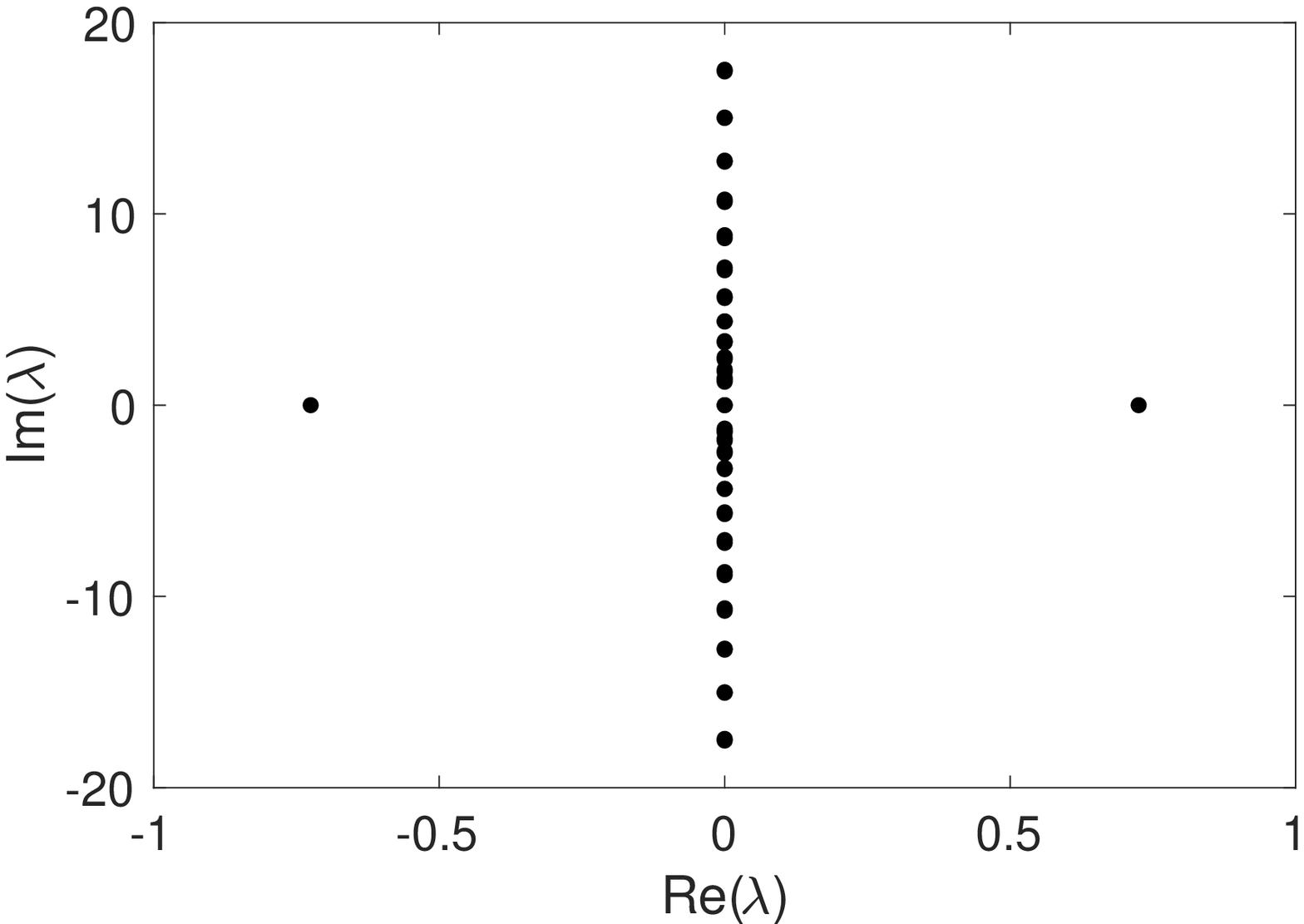}}\\
	\subfloat[]{\includegraphics[scale=0.4]{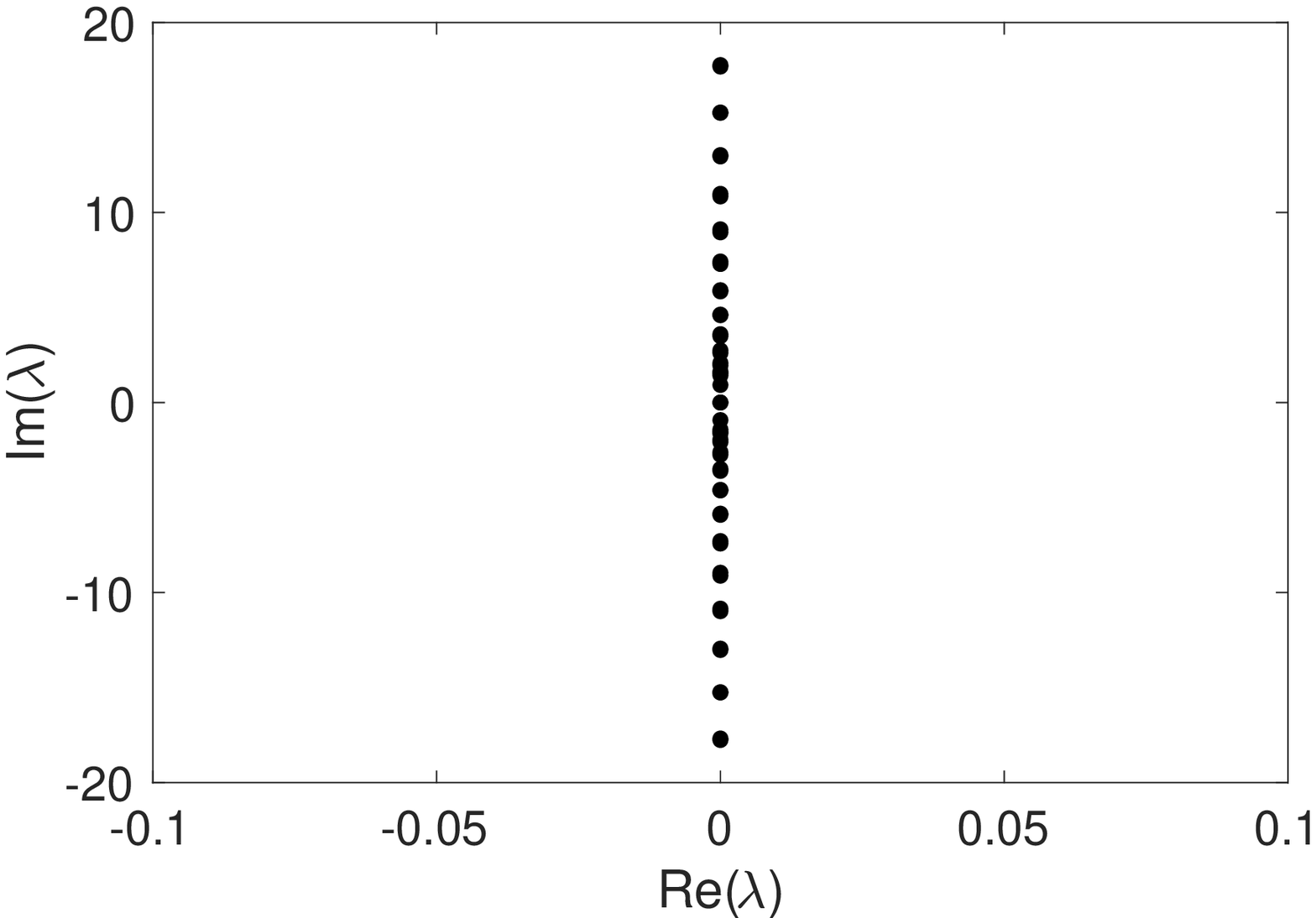}}
	\subfloat[]{\includegraphics[scale=0.4]{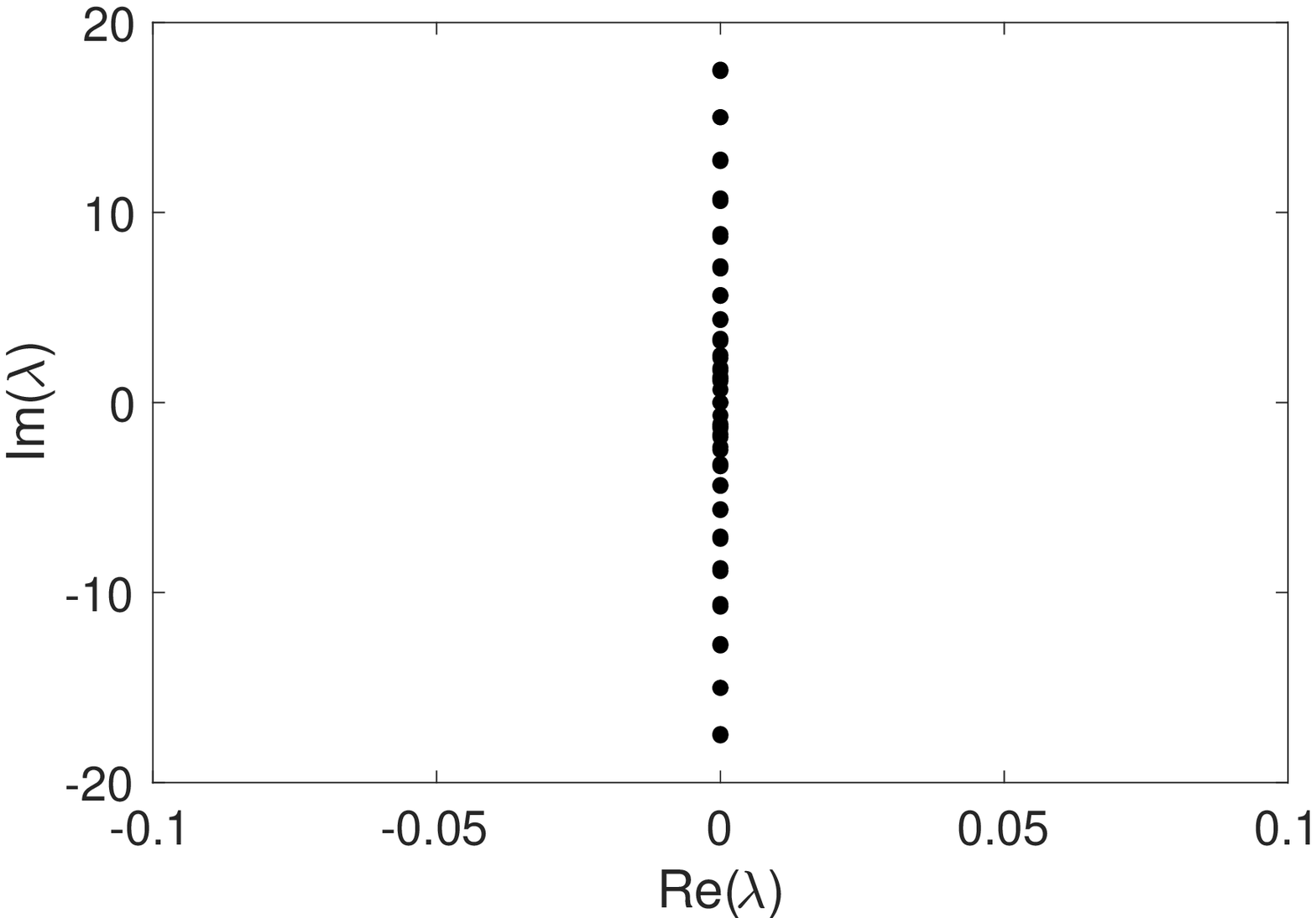}}
	\caption{The same as Fig.\ \ref{fig:sp_sol2}, but for the solitons in Fig.\ \ref{fig:sol2}. 
	}
	\label{fig:sp_sol3}
\end{figure}

\subsection{when $x>L$}
In this region the base point of the vector $\mathbf{b}$ will be carried under the flow to the origin. We will determine whether the landing curve $J^{-1}(W^s)$, intersecting the inner orbit, yields an additional vertical crossing or not. First, we look at the case $L_1 < \bar{L}_1$.Here the second defect the vector $\mathbf{b}$ still points to the transient orbit. 
%Comparing $\mathbf{b}$ with the vector that is taken to the curve $W^s$ which is tangent to $J^{-1}(W^s)$, 
We can see that $\mathbf{b}$ is lower than the vector that is taken to the curve $W^s$ which is tangent to $J^{-1}(W^s)$ (i.e. $\mathbf{b}$ has a larger angle). Comparing these two vectors, in this region, there will be no additional crossing to the vertical. 

Now, for the case $L_1 \geq \bar{L}_1$. If $L_2>\bar{L}_2$ at the second defect, the vector $\mathbf{b}$ is pointing out from the transient orbit relative to the vector that is tangent to $J^{-1}(W^s)$. So $\mathbf{b}$ has a smaller angle. After the jump, the flow pushes it across the vertical, so in this case $P \geq 2$. If $L_2 \leq \bar{L}_2$, $\mathbf{b}$ has a larger angle, and there are no additional crossings of the vertical.

To summarise, we have the following 
\begin{thm}
	Positive definite homoclinic solutions of \eqref{stateq}-\eqref{mm} with $L_1 < \bar{L}_1$ will have $P \leq 1$. If they have $L_1 \geq \bar{L}_1$, then there are two possible cases, i.e., either $L_2 < \bar{L}_2$ or $L_2 \geq \bar{L}_2$. The former case gives $P \leq 1$, while the latter yields $P \geq 2$. 
\end{thm}
Using Theorem \ref{th1}, the last case will give an unstable solution through a real eigenvalue. {Solutions in Figs.\ \ref{subfig:sym_L_1_e_1_w_0k6}, \ref{subfig:asym_L_1_e_1_w_1}, and \ref{subfig:S3_L_1_e_0k95_w_1k4} correspond to $L_1 < \bar{L_1}$. Solutions in Figs.\ \ref{subfig:S2_L_1_e_1_w_1}, \ref{subfig:S2_L_1_e_0k95_w_0k7} and \ref{subfig:S2_L_1_e_0k95_w_1k1} correspond to $L_1 > \bar{L_1}$, but $L_2<\bar{L_2}$. In those cases, we cannot determine their stability. Using numerics, our results in the next section show that they are stable. On the other hand, for the solutions in Fig.\ \ref{subfig:sym_L_1_e_1_w_1k2} and \ref{subfig:S1_L_1_e_0k95_w_1k2},  $L_1 > \bar{L_1}$ and at the same time $L_2>\bar{L_2}$. Hence, they are unstable. %For the case $\epsilon=0.95$, the solutions in Fig~ are stable because $L_1 > \bar{L_1}$ and $L_2 < \bar{L_2}$, the solution in Fig~ is unstable because $L_1 > \bar{L_1}$ and $L_2>\bar{L_2}$.}

%\textbf{Tolong dilengkapi dan dihubungkan dengan solusi yang diplot di Fig 2 dan 3. Mana-mana yang $L_1$ atau $L_2$, dst. Sehingga kesimpulannya mana yang stabil dan mana yang tidak dalam gambar 2 dan 3 itu.}

\section{Numerical results}
\label{num}

We solved Eqs.\ \eqref{stateq} and \eqref{mm} as well as Eq.\ \eqref{sys} numerically to study the localised standing waves and their stability. A central finite difference was used to approximate the Laplacian with a relatively fine discretisation. Here we  present the spectrum of the solutions from \ref{stab} obtained from solving the eigenvalue problem numerically. 

We plot the spectrum of solutions in Figs.\ \ref{fig:sol} and \ref{fig:sol2} in Figs.\ \ref{fig:sp_sol2} and \ref{fig:sp_sol3}, respectively. We confirm the result of Section \ref{stab} that solutions plotted in panel (b) of Figs.\ \ref{fig:sol} and \ref{fig:sol2} are unstable. The instability is due to the presence of a pair of real eigenvalues. 

\begin{figure}[htbp!]
	\centering
	{\includegraphics[scale=0.6]{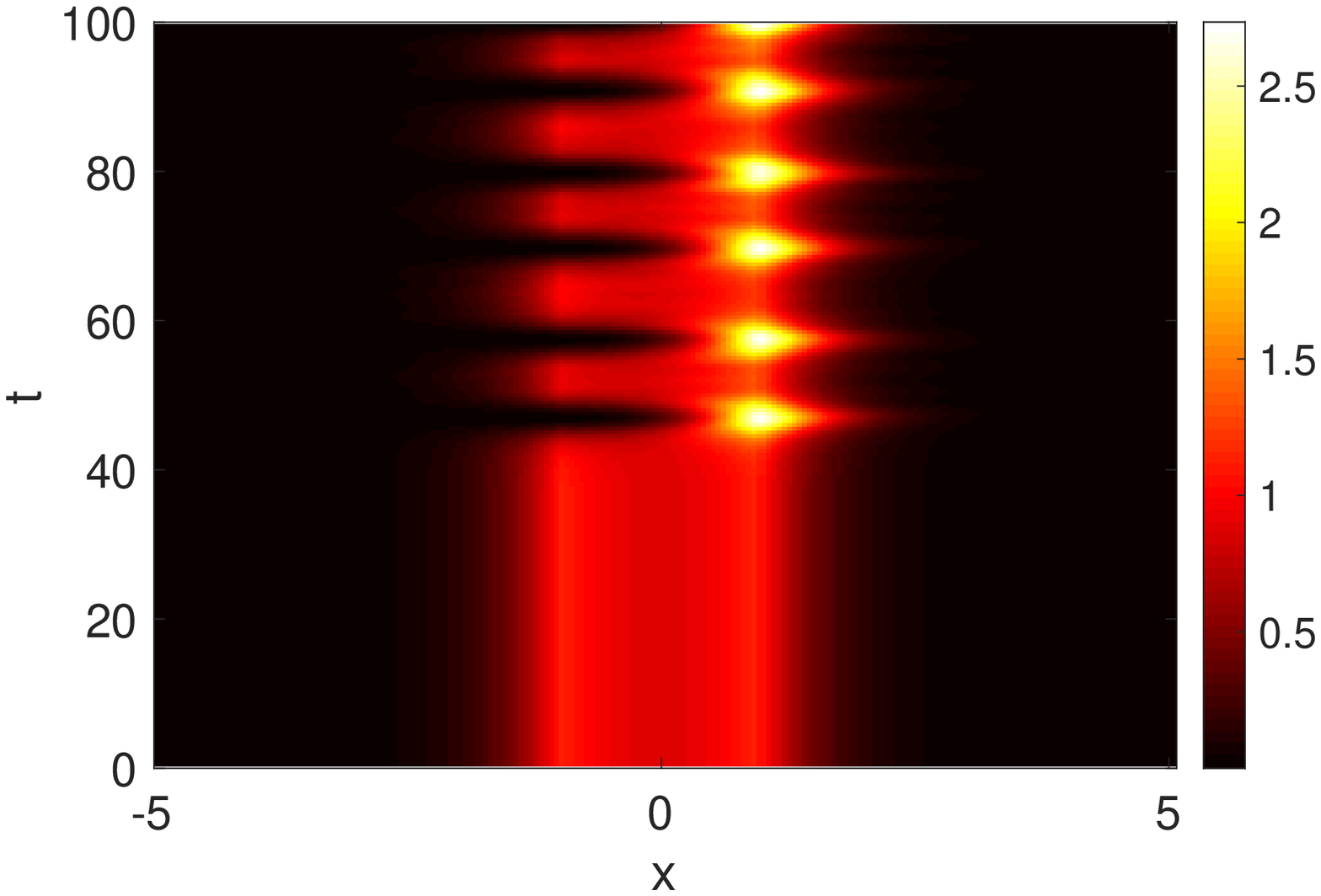}}
	\caption{Time dynamics of the unstable solution in Fig.\ \ref{subfig:S1_L_1_e_0k95_w_1k2}. Plotted is the squared magnitude $|\psi|^2$. Initially the standing wave is perturbed randomly.
	}
	\label{fig:dyn}
\end{figure}

When a solution is unstable, it is interesting to see its typical dynamics. To do so, we solve the governing equation \eqref{nls} numerically where the Dirac delta potential is incorporated through the boundaries. While the spatial discretisation is still the same as before, the time derivative is integrated using the classic fourth-order Runge-Kutta method. 

In Fig.\ \ref{fig:dyn} we plot time dynamic of the unstable solution shown in panel (b) of Fig.\ \ref{fig:sol2}. The time evolution is typical where the instability manifests in the form of periodic oscillations. The norm tends to be localised in one of the wells, which is one of the characteristics of the presence of symmetry breaking solutions \cite{kevrekidis2005spontaneous,marzuola2010long,susanto11,susanto12,theo}. 

\section{Conclusion}
\label{concl}

In this paper, we have considered broken symmetry breaking bifurcations in the NLS on the real line with an asymmetric double Dirac delta potential. By using a dynamical system approach, we presented the ground state solutions in the phase plane and their explicit expressions. We have shown that different from the symmetric case where the bifurcation is of a pitchfork type, when the potential is asymmetric, the bifurcation is of a saddle-centre type. The linear instability of the corresponding solutions has been derived as well using a geometrical approach developed by Jones \cite{jones1988instability}. Numerical computations have been presented illustrating the analytical results and simulations showing the typical dynamics of unstable solutions have also been discussed.

\section*{Acknowledgments}
R.R gratefully acknowledges financial support from Lembaga Pengelolaan Dana Pendidikan (Indonesia Endowment Fund for Education), Grant Ref.\ No: S-5405/LPDP.3/2015. \\

The authors contributed equally to the manuscript. \\

%\section*{References}
%
%\bibliography{mybibfile}

\begin{thebibliography}{10}
\bibliographystyle{amsalpha}

\bibitem{part}
TWB~Kibble.  \newblock Spontaneous symmetry breaking in gauge theories. \newblock {\em Philosophical Transactions of the Royal Society A: Mathematical, Physical and Engineering Sciences} 373(2032) : 20140033, 2015. 
\bibitem{albiez2005direct}
Michael Albiez, Rudolf Gati, Jonas F{\"o}lling, Stefan Hunsmann, Matteo
Cristiani, and Markus~K Oberthaler.
\newblock Direct observation of tunneling and nonlinear self-trapping in a
single bosonic Josephson junction.
\newblock {\em Physical Review Letters}, 95(1):010402, 2005.

\bibitem{zibold2010classical}
Tilman Zibold, Eike Nicklas, Christian Gross, and Markus~K Oberthaler.
\newblock Classical bifurcation at the transition from Rabi to Josephson
dynamics.
\newblock {\em Physical Review Letters}, 105(20):204101, 2010.

\bibitem{liu2014spontaneous}
Mingkai Liu, David~A Powell, Ilya~V Shadrivov, Mikhail Lapine, and Yuri~S
Kivshar.
\newblock Spontaneous chiral symmetry breaking in metamaterials.
\newblock {\em Nature Communications}, 5:4441, 2014.

\bibitem{green1990spontaneous}
C~Green, GB~Mindlin, EJ~D’Angelo, HG~Solari, and JR~Tredicce.
\newblock Spontaneous symmetry breaking in a laser: the experimental side.
\newblock {\em Physical Review Letters}, 65(25):3124, 1990.

\bibitem{kevrekidis2005spontaneous}
PG~Kevrekidis, Zhigang Chen, BA~Malomed, DJ~Frantzeskakis, and MI~Weinstein.
\newblock Spontaneous symmetry breaking in photonic lattices: Theory and
experiment.
\newblock {\em Physics Letters A}, 340(1-4):275--280, 2005.

\bibitem{sawai2000spontaneous}
Satoshi Sawai, Yasuo Maeda, and Yasuji Sawada.
\newblock Spontaneous symmetry breaking turing-type pattern formation in a
confined dictyostelium cell mass.
\newblock {\em Physical Review Letters}, 85(10):2212, 2000.

\bibitem{heil2001chaos}
Tilmann Heil, Ingo Fischer, Wolfgang Els{\"a}sser, Josep Mulet, and Claudio~R
Mirasso.
\newblock Chaos synchronization and spontaneous symmetry-breaking in
symmetrically delay-coupled semiconductor lasers.
\newblock {\em Physical Review Letters}, 86(5):795, 2001.

\bibitem{hamel2015spontaneous}
Philippe Hamel, Samir Haddadi, Fabrice Raineri, Paul Monnier, Gregoire
Beaudoin, Isabelle Sagnes, Ariel Levenson, and Alejandro~M Yacomotti.
\newblock Spontaneous mirror-symmetry breaking in coupled photonic-crystal
nanolasers.
\newblock {\em Nature Photonics}, 9(5):311, 2015.

\bibitem{davies1979symmetry}
EB~Davies.
\newblock Symmetry breaking for a non-linear Schr{\"o}dinger equation.
\newblock {\em Communications in Mathematical Physics}, 64(3):191--210, 1979.

\bibitem{mahmud2002bose}
KW~Mahmud, JN~Kutz, and WP~Reinhardt.
\newblock Bose-Einstein condensates in a one-dimensional double square well:
Analytical solutions of the nonlinear Schr{\"o}dinger equation.
\newblock {\em Physical Review A}, 66(6):063607, 2002.

\bibitem{marzuola2010long}
Jeremy~L Marzuola and Michael~I Weinstein.
\newblock Long time dynamics near the symmetry breaking bifurcation for
nonlinear Schr{\"o}dinger/Gross-Pitaevskii equations.
\newblock {\em Discrete \& Continuous Dynamical Systems-A}, 28(4):1505--1554,
2010.

\bibitem{jackson2004geometric}
Russell~K Jackson and Michael~I Weinstein.
\newblock Geometric analysis of bifurcation and symmetry breaking in a
Gross—Pitaevskii equation.
\newblock {\em Journal of Statistical Physics}, 116(1-4):881--905, 2004.

\bibitem{susanto11}
H~Susanto, J~Cuevas, and P~Kr\"uger. 
\newblock Josephson tunnelling of dark solitons in a double-well potential. 
\newblock \emph{Journal of Physics B: Atomic, Molecular and Optical Physics} 44(9):095003, 2011.

\bibitem{susanto12}
H~Susanto and J~Cuevas. 
\newblock Josephson tunneling of excited states in a double-well potential. \newblock In \emph{Spontaneous Symmetry Breaking, Self-Trapping, and Josephson Oscillations}. Springer, Berlin, Heidelberg, 2012. 583-599.

\bibitem{theo}
G~Theocharis, PG~Kevrekidis, DJ~Frantzeskakis, and P~Schmelcher. 
\newblock Symmetry breaking in symmetric and asymmetric double-well potentials. 
\newblock \emph{Physical Review E} 74(5): 056608, 2006.

\bibitem{jones1988instability}
Christopher~KRT Jones.
\newblock Instability of standing waves for non-linear Schr{\"o}dinger-type
equations.
\newblock {\em Ergodic Theory and Dynamical Systems}, 8(8*):119--138, 1988.

\bibitem{mara1} 
R~Marangell, CKRT~Jones, and H~Susanto.
\newblock Localized standing waves in inhomogeneous Schr{\"o}dinger equations.
\newblock {\em Nonlinearity}, 23(9):2059, 2010.
\bibitem{mara2}
R~Marangell, H~Susanto, and CKRT~Jones. \newblock Unstable gap solitons in inhomogeneous nonlinear Schr\"odinger equations. \newblock \emph{Journal of Differential Equations} 253(4): 1191-1205,  2012.


\bibitem{nickalls1993new}
Richard~WD Nickalls.
\newblock A new approach to solving the cubic: Cardan’s solution revealed.
\newblock {\em The Mathematical Gazette}, 77(480):354--359, 1993.

\bibitem{landa2001regular}
Polina~S Landa.
\newblock {\em Regular and Chaotic Oscillations}.
\newblock Springer Science \& Business Media, 2001.

\end{thebibliography}

%\section*{References}

\end{document}